\documentclass[12pt]{article}
\usepackage{theorem}
\usepackage{amssymb}
\usepackage{graphicx}
\usepackage[mathscr]{eucal}
\newtheorem{prop}{}[section]

{\theorembodyfont{\upshape} \newtheorem{rema}[prop]{}}
\newcommand{\boma}[1]{{\mbox{\boldmath $#1$} }}

\begin{document}
\newcommand{\uper}[1]{\stackrel{\barray{c} {~} \\ \mbox{\footnotesize{#1}}\farray}{\longrightarrow} }
\newcommand{\nop}[1]{ \|#1\|_{\piu} }
\newcommand{\no}[1]{ \|#1\| }
\newcommand{\nom}[1]{ \|#1\|_{\meno} }
\newcommand{\UU}[1]{e^{#1 \AA}}
\newcommand{\UD}[1]{e^{#1 \Delta}}
\newcommand{\bb}[1]{\mathbb{{#1}}}
\newcommand{\HO}[1]{\bb{H}^{{#1}}}
\newcommand{\Hz}[1]{\bb{H}^{{#1}}_{\zz}}
\newcommand{\Hs}[1]{\bb{H}^{{#1}}_{\ss}}
\newcommand{\Hg}[1]{\bb{H}^{{#1}}_{\gg}}
\newcommand{\HM}[1]{\bb{H}^{{#1}}_{\so}}
\def\tvainf{\vspace{-0.4cm} \barray{ccc} \vspace{-0,1cm}{~}
\\ \vspace{-0.2cm} \longrightarrow \\ \vspace{-0.2cm} \scriptstyle{T \vain + \infty} \farray}
\def\comple{\scriptscriptstyle{\complessi}}
\def\nume{0.407}
\def\numerob{0.00724}
\def\deln{7/10}
\def\delnn{\dd{7 \over 10}}
\def\e{c}
\def\p{p}
\def\z{z}
\def\symd{{\mathfrak S}_d}
\def\del{\omega}
\def\Del{\delta}
\def\Di{\Delta}
\def\Ss{{\mathscr{S}}}
\def\mmu{\hat{\mu}}
\def\rot{\mbox{rot}\,}
\def\curl{\mbox{curl}\,}
\def\Mm{\mathscr M}
\def\XS{\boma{x}}
\def\TS{\boma{t}}
\def\Lam{\boma{\eta}}
\def\DS{\boma{\rho}}
\def\KS{\boma{k}}
\def\LS{\boma{\lambda}}
\def\PR{\boma{p}}
\def\VS{\boma{v}}
\def\ski{\! \! \! \! \! \! \! \! \! \! \! \! \! \!}
\def\h{L}
\def\EM{M}
\def\EMP{M'}
\def\R{R}
\def\Rr{\Upsilon}
\def\E{E}
\def\FFf{\mathscr{F}}
\def\A{F}
\def\Xim{\Xi_{\meno}}
\def\Ximn{\Xi_{n-1}}
\def\lan{\lambda}
\def\om{\omega}
\def\Om{\Omega}
\def\Sim{\Sigm}
\def\Sip{\Delta \Sigm}
\def\Sigm{{\mathscr{S}}}
\def\Ki{{\mathscr{K}}}
\def\Hi{{\mathscr{H}}}
\def\zz{{\scriptscriptstyle{0}}}
\def\ss{{\scriptscriptstyle{\Sigma}}}
\def\gg{{\scriptscriptstyle{\Gamma}}}
\def\so{\ss \zz}
\def\Dz{\bb{\DD}'_{\zz}}
\def\Ds{\bb{\DD}'_{\ss}}
\def\Dsz{\bb{\DD}'_{\so}}
\def\Dg{\bb{\DD}'_{\gg}}
\def\Ls{\bb{L}^2_{\ss}}
\def\Lg{\bb{L}^2_{\gg}}
\def\bF{{\bb{V}}}
\def\Fz{\bF_{\zz}}
\def\Fs{\bF_\ss}
\def\Fg{\bF_\gg}
\def\Pre{P}
\def\UUU{{\mathcal U}}
\def\fiapp{\phi}
\def\PU{P1}
\def\PD{P2}
\def\PT{P3}
\def\PQ{P4}
\def\PC{P5}
\def\PS{P6}
\def\Q{P6}
\def\X{Q2}
\def\Xp{Q3}
\def\Vi{V}
\def\bVi{\bb{V}}
\def\K{V}
\def\Ks{\bb{\K}_\ss}
\def\Kz{\bb{\K}_0}
\def\KM{\bb{\K}_{\, \so}}
\def\HGG{\bb{H}^\G}
\def\HG{\bb{H}^\G_{\so}}
\def\EG{{\mathfrak{P}}^{\G}}
\def\G{G}
\def\de{\delta}
\def\esp{\sigma}
\def\dd{\displaystyle}
\def\LP{\mathfrak{L}}
\def\dive{\mbox{div}}
\def\la{\langle}
\def\ra{\rangle}
\def\um{u_{\meno}}
\def\uv{\mu_{\meno}}
\def\Fp{ {\textbf F_{\piu}} }
\def\Ff{ {\textbf F} }
\def\Fm{ {\textbf F_{\meno}} }
\def\piu{\scriptscriptstyle{+}}
\def\meno{\scriptscriptstyle{-}}
\def\omeno{\scriptscriptstyle{\ominus}}
\def\Tt{ {\mathscr T} }
\def\Xx{ {\textbf X} }
\def\Yy{ {\textbf Y} }
\def\Ee{ {\textbf E} }
\def\VP{{\mbox{\tt VP}}}
\def\CP{{\mbox{\tt CP}}}
\def\cp{$\CP(f_0, t_0)\,$}
\def\cop{$\CP(f_0)\,$}
\def\copn{$\CP_n(f_0)\,$}
\def\vp{$\VP(f_0, t_0)\,$}
\def\vop{$\VP(f_0)\,$}
\def\vopn{$\VP_n(f_0)\,$}
\def\vopdue{$\VP_2(f_0)\,$}
\def\leqs{\leqslant}
\def\geqs{\geqslant}
\def\mat{{\frak g}}
\def\tG{t_{\scriptscriptstyle{G}}}
\def\tN{t_{\scriptscriptstyle{N}}}
\def\TK{t_{\scriptscriptstyle{K}}}
\def\CK{C_{\scriptscriptstyle{K}}}
\def\CN{C_{\scriptscriptstyle{N}}}
\def\CG{C_{\scriptscriptstyle{G}}}
\def\CCG{{\mathscr{C}}_{\scriptscriptstyle{G}}}
\def\tf{{\tt f}}
\def\ti{{\tt t}}
\def\ta{{\tt a}}
\def\tc{{\tt c}}
\def\tF{{\tt R}}
\def\C{{\mathscr C}}
\def\P{{\mathscr P}}
\def\V{{\mathscr V}}
\def\TI{\tilde{I}}
\def\TJ{\tilde{J}}
\def\Lin{\mbox{Lin}}
\def\Hinfc{ H^{\infty}(\reali^d, \complessi) }
\def\Hnc{ H^{n}(\reali^d, \complessi) }
\def\Hmc{ H^{m}(\reali^d, \complessi) }
\def\Hac{ H^{a}(\reali^d, \complessi) }
\def\Dc{\DD(\reali^d, \complessi)}
\def\Dpc{\DD'(\reali^d, \complessi)}
\def\Sc{\SS(\reali^d, \complessi)}
\def\Spc{\SS'(\reali^d, \complessi)}
\def\Ldc{L^{2}(\reali^d, \complessi)}
\def\Lpc{L^{p}(\reali^d, \complessi)}
\def\Lqc{L^{q}(\reali^d, \complessi)}
\def\Lrc{L^{r}(\reali^d, \complessi)}
\def\Hinfr{ H^{\infty}(\reali^d, \reali) }
\def\Hnr{ H^{n}(\reali^d, \reali) }
\def\Hmr{ H^{m}(\reali^d, \reali) }
\def\Har{ H^{a}(\reali^d, \reali) }
\def\Dr{\DD(\reali^d, \reali)}
\def\Dpr{\DD'(\reali^d, \reali)}
\def\Sr{\SS(\reali^d, \reali)}
\def\Spr{\SS'(\reali^d, \reali)}
\def\Ldr{L^{2}(\reali^d, \reali)}
\def\Hinfk{ H^{\infty}(\reali^d, \KKK) }
\def\Hnk{ H^{n}(\reali^d, \KKK) }
\def\Hmk{ H^{m}(\reali^d, \KKK) }
\def\Hak{ H^{a}(\reali^d, \KKK) }
\def\Dk{\DD(\reali^d, \KKK)}
\def\Dpk{\DD'(\reali^d, \KKK)}
\def\Sk{\SS(\reali^d, \KKK)}
\def\Spk{\SS'(\reali^d, \KKK)}
\def\Ldk{L^{2}(\reali^d, \KKK)}
\def\Knb{K^{best}_n}
\def\sc{\cdot}
\def\k{\mbox{{\tt k}}}
\def\x{\mbox{{\tt x}}}
\def\g{ {\textbf g} }
\def\QQQ{ {\textbf Q} }
\def\AAA{ {\textbf A} }
\def\gr{\mbox{gr}}
\def\sgr{\mbox{sgr}}
\def\loc{\mbox{loc}}
\def\PZ{{\Lambda}}
\def\PZAL{\mbox{P}^{0}_\alpha}
\def\epsilona{\epsilon^{\scriptscriptstyle{<}}}
\def\epsilonb{\epsilon^{\scriptscriptstyle{>}}}
\def\lgraffa{ \mbox{\Large $\{$ } \hskip -0.2cm}
\def\rgraffa{ \mbox{\Large $\}$ } }
\def\restriction{\upharpoonright}
\def\M{{\scriptscriptstyle{M}}}
\def\m{m}
\def\Fre{Fr\'echet~}
\def\I{{\mathcal N}}
\def\ap{{\scriptscriptstyle{ap}}}
\def\fiap{\varphi_{\ap}}
\def\dfiap{{\dot \varphi}_{\ap}}
\def\DDD{ {\mathfrak D} }
\def\BBB{ {\textbf B} }
\def\EEE{ {\textbf E} }
\def\GGG{ {\textbf G} }
\def\TTT{ {\textbf T} }
\def\KKK{ {\textbf K} }
\def\HHH{ {\textbf K} }
\def\FFi{ {\bf \Phi} }
\def\GGam{ {\bf \Gamma} }
\def\sc{ {\scriptstyle{\bullet} }}
\def\a{a}
\def\ep{\epsilon}
\def\c{\kappa}
\def\parn{\par\noindent}
\def\teta{M}
\def\elle{L}
\def\ro{\rho}
\def\al{\alpha}
\def\si{\sigma}
\def\be{\beta}
\def\ga{\gamma}
\def\te{\vartheta}
\def\ch{\chi}
\def\et{\eta}
\def\complessi{{\bf C}}
\def\len{{\bf L}}
\def\reali{{\bf R}}
\def\interi{{\bf Z}}
\def\Z{{\bf Z}}
\def\naturali{{\bf N}}
\def\To{ {\bf T} }
\def\Td{ {\To}^d }
\def\Tt{ {\To}^3 }
\def\Zd{ \interi^d }
\def\Zt{ \interi^3 }
\def\Zet{{\mathscr{Z}}}
\def\Ze{\Zet^d}
\def\T1{{\textbf To}^{1}}
\def\es{s}
\def\ee{{E}}
\def\FF{\mathcal F}
\def\FFu{ {\textbf F_{1}} }
\def\FFd{ {\textbf F_{2}} }
\def\GG{{\mathcal G} }
\def\EE{{\mathcal E}}
\def\KK{{\mathcal K}}
\def\PP{{\mathcal P}}
\def\PPP{{\mathscr P}}
\def\PN{{\mathcal P}}
\def\PPN{{\mathscr P}}
\def\QQ{{\mathcal Q}}
\def\J{J}
\def\Np{{\hat{N}}}
\def\Lp{{\hat{L}}}
\def\Jp{{\hat{J}}}
\def\Pp{{\hat{P}}}
\def\Pip{{\hat{\Pi}}}
\def\Vp{{\hat{V}}}
\def\Ep{{\hat{E}}}
\def\Gp{{\hat{G}}}
\def\Kp{{\hat{K}}}
\def\Ip{{\hat{I}}}
\def\Tp{{\hat{T}}}
\def\Mp{{\hat{M}}}
\def\La{\Lambda}
\def\Ga{\Gamma}
\def\Si{\Sigma}
\def\Upsi{\Upsilon}
\def\Gam{\Gamma}
\def\Gag{{\check{\Gamma}}}
\def\Lap{{\hat{\Lambda}}}
\def\Upsig{{\check{\Upsilon}}}
\def\Kg{{\check{K}}}
\def\ellp{{\hat{\ell}}}
\def\j{j}
\def\jp{{\hat{j}}}
\def\BB{{\mathcal B}}
\def\LL{{\mathcal L}}
\def\MM{{\mathcal U}}
\def\SS{{\mathcal S}}
\def\DD{D}
\def\Dd{{\mathcal D}}
\def\VV{{\mathcal V}}
\def\WW{{\mathcal W}}
\def\OO{{\mathcal O}}
\def\RR{{\mathcal R}}
\def\TT{{\mathcal T}}
\def\AA{{\mathcal A}}
\def\CC{{\mathcal C}}
\def\JJ{{\mathcal J}}
\def\NN{{\mathcal N}}
\def\HH{{\mathcal H}}
\def\XX{{\mathcal X}}
\def\XXX{{\mathscr X}}
\def\YY{{\mathcal Y}}
\def\ZZ{{\mathcal Z}}
\def\CC{{\mathcal C}}
\def\cir{{\scriptscriptstyle \circ}}
\def\circa{\thickapprox}
\def\vain{\rightarrow}
\def\parn{\par \noindent}
\def\salto{\vskip 0.2truecm \noindent}
\def\spazio{\vskip 0.5truecm \noindent}
\def\vs1{\vskip 1cm \noindent}
\def\fine{\hfill $\square$ \vskip 0.2cm \noindent}
\def\ffine{\hfill $\lozenge$ \vskip 0.2cm \noindent}
\newcommand{\rref}[1]{(\ref{#1})}
\def\beq{\begin{equation}}
\def\feq{\end{equation}}
\def\beqq{\begin{eqnarray}}
\def\feqq{\end{eqnarray}}
\def\barray{\begin{array}}
\def\farray{\end{array}}
\makeatletter \@addtoreset{equation}{section}
\renewcommand{\theequation}{\thesection.\arabic{equation}}
\makeatother
\begin{titlepage}
{~}
\vspace{-2cm}
\begin{center}
{\huge An $H^1$ setting for the Navier-Stokes equations: quantitative estimates}
\end{center}
\vspace{0.5truecm}
\begin{center}
{\large
Carlo Morosi${}^1$, Livio Pizzocchero${}^2$} \\
\vspace{0.5truecm} ${}^1$ Dipartimento di Matematica, Politecnico di Milano,
\\ P.za L. da Vinci 32, I-20133 Milano, Italy \\
e--mail: carlo.morosi@polimi.it \\
${}^2$ Dipartimento di Matematica, Universit\`a di Milano\\
Via C. Saldini 50, I-20133 Milano, Italy\\
and Istituto Nazionale di Fisica Nucleare, Sezione di Milano, Italy \\
e--mail: livio.pizzocchero@unimi.it
\end{center}
\begin{abstract}
We consider the incompressible Navier-Stokes (NS) equations on a torus, in the
setting of the spaces $L^2$ and $H^1$; our approach is
based on a general framework for semi- or quasi-linear parabolic equations
proposed in the previous work \cite{due}. We present some estimates
on the linear semigroup generated by the Laplacian and on
the quadratic NS nonlinearity; these are fully quantitative, i.e., all the
constants appearing therein are given explicitly. \parn
As an application we show that, on a three dimensional torus $\Tt$, the (mild) solution of
the NS Cauchy problem is global for each $H^1$ initial datum $u_0$ with zero mean, such that
$|| \curl u_0 ||_{L^2} \leqs \nume$; this improves the bound for global existence
$\| \curl u_0 ||_{L^2} \leqs \numerob$,
derived recently by Robinson and Sadowski \cite{Rob}. We announce some future applications,
based again on the $H^1$ framework and on the general scheme of \cite{due}.
\end{abstract}
\vspace{0.2cm} \noindent
\textbf{Keywords:} Differential equations, theoretical approximation, Navier-Stokes \hfill \parn equations.
\par \vspace{0.05truecm} \noindent \textbf{AMS 2000 Subject classifications:} 35Q30, 76D03, 76D05.
\end{titlepage}
\section{Introduction}
As well known, the great open problem about the incompressible Navier-Stokes (NS) equations is to prove
global existence of the solutions in three space dimensions, for all
sufficiently smooth initial data. \parn
A much more modest, but realistic research program about
these equations is the following. \parn
(i) Deriving estimates for the linear semigroup generated by the Laplacian
and for the NS quadratic nonlinearity, within
an appropriate functional setting (say, the Sobolev spaces
on $\reali^3$ or on the torus $\Tt$). This should be done paying attention to the
strictly quantitative aspects, and having as a goal the best achievable accuracy.
\parn
(ii) Explicitating the estimates
for the NS solutions derivable from the previous framework.
\par \noindent
Some typical issues to be treated as parts of item (ii), always using
the information from (i), are the following ones: \parn
(ii.a) Deriving quantitative lower bounds on the interval of existence of the solution,
for a given initial datum. \parn
(ii.b) Producing a sufficient condition for global existence of the solution, when the
norm of the initial datum is below an \textsl{explicitly given} upper bound. \parn
(ii.c) Devising \textsl{a posteriori} tests
on the approximate solutions (e.g., of Galerkin type) to get
quantitative  estimates on exact solutions.
\vskip 0.2cm \noindent
A number of papers published in recent times are more or less related to the above research
program. We mention, in particular:
the works by Chernyshenko, Constantin, Robinson,
Titi \cite{Che}, Dashti and Robinson \cite{DR} on the approximate solutions
of the NS equations, and on the \textsl{a posteriori} estimates derivable from them; the paper
by Robinson and Sadowski \cite{Rob} on the conditions for the existence of (exact) global
solutions, containing very interesting considerations on the computational times required to check them via
suitable approximate solutions; our papers \cite{uno} \cite{due} on approximate solutions
and \textsl{a posteriori} estimates for general semi- or quasi-linear evolution equations,
including (in the second reference) some applications to NS equations on a torus. Let
us also mention that a similar attitude towards approximate solutions and \textsl{a posteriori} bounds
has been developed by Machiels, Peraire and Patera \cite{Pat}, Veroy and Patera \cite{Pat2}
in connection with the NS equations
(or some related space discretizations), dealing mainly with two-dimensional steady problems.
\vskip 0.2cm \noindent
The present work outlines a framework to treat problems (ii.a)(ii.b)(ii.c) for the incompressible
NS equations on the torus $\Tt$. The main application presented here is a fully quantitative bound
on the initial datum yielding global existence, in the spirit of (ii.b);
a few words will be spent on other issues within the scheme of (ii), to be treated elsewhere. \parn
Our present constructions are related to the setup of \cite{due}, where
the basic functional space for the three dimensional NS equations
was the Sobolev space $\HM{n}(\Tt) \subset \Hz{n}(\Tt)$, with $n > 3/2$ (the subscripts $\Sigma$, $0$
indicating the conditions of zero divergence and zero mean, respectively). \parn
On the contrary, here the solutions of the NS equations
take values in $\HM{1}(\Tt) \subset \Hz{1}(\Tt)$. $\Hz{1}(\Tt)$
is equipped with the norm $v \mapsto \| v\|_1 := \| \sqrt{-\Delta} \, v \|_{L^2}$  ($\Delta$ the
Laplacian); the latter is restricted to $\HM{1}(\Tt)$, where it is found to coincide with the $L^2$ norm
of the vorticity:
\beq \| v \|_{1} = \| \curl v \|_{L^2}~. \feq
Throughout the paper, the incompressible NS equations for the velocity field $u$ are written as
\beq {d u \over d t}  = \Delta u - \LP(u \sc \partial u)~, \feq
with $\LP$ indicating the Leray projection on the divergence free vector fields.
As in our previous work, we emphasize the role played by the constants in certain basic estimates
on the semigroup $(e^{t \Delta})$ and on the NS bilinear map $(v,w) \mapsto \LP(v \sc \partial w)$. In particular,
for $1/2 < \del < 1$ we consider the negative order Sobolev space
$\HM{-\del}(\Tt)$ with the norm $\| \cdot \|_{-\del} := \| \sqrt{\Delta}^{\, -\del} \cdot \,\|_{L^2}$, 
and give explicitly a constant $K_{\del}$ such that, for all $v,w \in \HM{1}(\Tt)$, 
\beq \| \LP(v \sc \partial w) \|_{-\del} \leqs K_{\del} \| v \|_1
\| w \|_1~; \label{kom} \feq
for example, if $\del = \deln$ the above inequality
is fulfilled with $K_{\deln} = 0.361$. On the grounds of our previous
investigations on multiplication in Sobolev spaces \cite{mult}, the constant
$K_{\om}$ determined for any $\om$ by our method is probably close to the best
constant fulfilling \rref{kom}. \parn
The evaluation of $K_{\om}$ can be combined with some strictly quantitative
estimates on $(e^{t \Delta})$, and with the general setting of
\cite{uno}, to obtain a number of results on the exact NS solutions. In particular, here
we show that the incompressible NS equations have a global
solution $u$ for any initial datum $u_0 \in \HM{1}(\Tt)$ such that
\beq \| u_0 \|_1 \leqs \nume~; \label{conume} \feq
this improves the quantitative condition for global existence presented in \cite{Rob} (page 43)
that, in our notations, can be written as $\| u_0 \|_1 \leqs \numerob$. In connection with
these results, let us mention that the idea of an explicit norm bound, sufficient for global existence,
dates back to the seminal paper by Kato and Fujita \cite{Kat}; here,
an approach slightly different from ours was outlined but not developed in a fully
quantitative way, since the authors did not compute all the necessary
constants (nor indicated how to do the missing calculations).
\vskip 0.2cm \noindent
Let us outline the organization of the present paper. In Section \ref{abf}, we review
some part of the abstract framework of \cite{due} for semi- or quasi-linear evolution equations;
we consider in particular the case of a quadratic nonlinearity, suitable for application
to the NS equations.
\parn
In Section \ref{inns}, we fix our standards about Sobolev spaces
on a torus $\Td$ (for the moment, of any dimension $d \geqs 2$), and indicate the connections
of this formalism with the NS equations (of course a precise definition of these
standards, and especially of the Sobolev norms, is required to make meaningful the quantitative estimates mentioned before,
to have comparisons with \cite{Rob} and other papers, etc.). \parn
In Section \ref{setting} we discuss the semigroup
$(t, v) \mapsto e^{t \Delta} v$ as a map $[0,+\infty)
\times \HM{1}(\Td) \vain \HM{1}(\Td)$, or as a map
$(0,+\infty) \times \HM{-\del}(\Td)$ $ \vain \HM{1}(\Td)$, for
$d \geqs 2$ and $0 < \del < 1$.
\parn
In Section \ref{sequad}, after some
preliminaries we discuss the bilinear map
$\HM{1}(\Td) \times \HM{1}(\Td)$ $\vain \HM{-\del}(\Td)$,
$(v, w) \mapsto \LP(v \sc \partial w)$ for $d \geqs 2$ and $\del > d/2-1$;
we infer an inequality of the type \rref{kom} and show how
to compute the related constant $K_{\om}$ (in fact, also dependent
on the dimension $d$). \parn
In Section \ref{toge}, we put together the results of Sections
\ref{setting}-\ref{sequad} on $(e^{t \Delta})$ and on the NS
bilinear map, which is possible if ($d= 2$, $0 < \del < 1$ or)
$d =3$, $1/2 < \del < 1$; this gives our final ``$H^1$ framework''
for the NS equations on $\Tt$. As an output of this framework,
choosing $\del = \deln$
we obtain the condition \rref{conume} of global existence in $\HM{1}(\Tt)$,
together with a quantitative estimate on the exponential
decay rate of the global solution $u$.
\parn
In Section \ref{develop}
we present a method to get
estimates on the (exact) NS solution analyzing \textsl{a posteriori} any approximate
solution, in the $H^1$ framework on $\Tt$. Here we specialize
the general method of the ''control inequality'' introduced
in \cite{uno} \cite{due}; this was already applied in \cite{due}
to the NS equations on $\HM{n}(\Tt)$, $n > 3/2$, for
example to get estimates from the Galerkin approximate solutions.
Some applications of the control inequality to the approximate
NS solutions in $\HM{1}(\Tt)$ will be presented elsewhere. \parn
The paper is completed by some Appendices, concerning more technical
issues. In Appendix \ref{appesp}, we prove some estimates
on the semigroup $(e^{t \Delta})$ employed in Sections \ref{inns}-\ref{setting}.
In Appendix \ref{appeco} we illustrate some general facts on (discrete) convolutions
of unimodal functions; these results are employed in the subsequent Appendix
\ref{appeka}, where we study a convolution whose maximum gives the constant
$K_{\del}$ of Eq. \rref{kom}. In Appendix \ref{appesup} we give some details
related to the computation of $K_{\del}$ in dimension $d=3$, for $\del = \deln$.
\section{An abstract framework for evolution equations with a quadratic nonlinearity}
\label{abf}
In \cite{due}, we considered a framework for evolutionary
problems (mainly, of parabolic type) with a nonlinear part of a fairly
general kind; in the cited paper, the framework
was subsequently specialized to the case of a quadratic nonlinearity. For our present purposes, it
suffices to review the case of a quadratic, time independent nonlinearity. \parn
\textbf{The framework.} Throughout this section, we consider a
set
\beq (\Fp, \Ff, \Fm, \AA, \PPP) \feq
with the following properties.
\vskip 0.1cm\noindent
(\PU) $\Fp$, $\Ff$ and $\Fm$ are Banach spaces with norms $\nop{~}$, $\no{~}$ and $\nom{~}$, such that
\beq \Fp \hookrightarrow \Ff \hookrightarrow \Fm \feq
(the symbol $\hookrightarrow$ indicating that one space is a dense linear subspace
of the other, and that the natural inclusion is continuous).
Elements of these spaces are generically denoted with $v, w,...$~. \parn
(\PD) $\AA$ is a linear operator such that
\beq \AA : \Fp \vain \Fm~, \qquad v \mapsto \AA v~.  \feq
Viewing $\Fp$ as a subspace of $\Fm$, the norm $\nop{~}$ is equivalent to the
graph norm $v \in \Fp \mapsto \nom{v} + \nom{A v}$. \vskip 0.1cm\noindent
(\PT) Viewing $\AA$ as a densely defined linear operator in $\Fm$,
it is assumed that $\AA$ generates a strongly continuous semigroup $(e^{t \AA})_{t \in [0,+\infty)}$
on $\Fm$ (of course, from the standard theory of linear semigroups, we have
$e^{t \AA}(\Fp) \subset \Fp$ for all $t \geqs 0$). \vskip 0.1cm\noindent
(\PQ) One has
\beq e^{t \AA}(\Ff) \subset \Ff \qquad \mbox{for $t \in [0, +\infty)$}~; \label{reguf} \feq
the function $(v, t) \mapsto e^{t \AA} v$ is continuous from $\Ff \times [0,+\infty)$ to $\Ff$, yielding
a strongly continuous semigroup on $\Ff$ as well.  Furthermore, there is a constant $B > 0$ such that
\beq \| e^{t \AA} v \| \leqs e^{- B t} \| v \| \qquad \mbox{for all $v \in \Ff$}~. \label{eqbi} \feq
(\PC) One has
\beq e^{t \AA}(\Fm) \subset \Ff \qquad \mbox{for $t \in(0, +\infty)$}~; \label{regum} \feq
the function $(v, t) \mapsto e^{t \AA} v$
is continuous from $\Fm \times (0,+\infty)$ to $\Ff$
(in a few words: for all $t > 0$, $e^{t \AA}$ regularizes the vectors of $\Fm$, sending them
into $\Ff$ continuously). Furthermore, there is a function
$\mu \in C((0,+\infty), (0,+\infty))$ such that
\beq \no{e^{t \AA} v} \leqs \mu(t) \, e^{- B t} \nom{v} \qquad \mbox{for $t > 0$, $v \in \Fm$, ~$B$ as in \rref{eqbi}}~.
\label{eqv} \feq
The function $\mu$ behaves like an \textsl{integrable} power of $1/t$ for $t \vain 0$, i.e.,
\beq \mu(t) = O({1 / t^{\varepsilon}}) \qquad \mbox{for $t \vain  0^{+}$}~,~~\varepsilon \in [0,1)~.
\label{integrab} \feq
Finally, there is a constant $N \in (0,+\infty)$ such that
\beq \int_{0}^t d s \, e^{-B s} \, \mu(t-s) \leqs N \qquad \mbox{for $t \in [0,+\infty)$}~.\label{defnn} \feq
\vskip 0.1cm\noindent
(\Q) $\PPP$ is a bilinear map such that
\beq \PPP : \Ff \times \Ff \vain \Fm~, \qquad (v,w) \mapsto \PPP(v,w)~; \feq
we assume continuity of $\PPP$, which is equivalent to the
existence of a constant $K \in [0,+\infty)$ such that, for all $v, w \in \Ff$,
\beq \nom{\PPP(v,w)} \leqs K \| v \| \| w \|~. \label{defk} \feq
\vskip 0.2cm \noindent
\textbf{The initial value problem; some results of uniqueness and existence.}
The initial value problem with initial datum $u_0 \in \Ff$ is the following one:
\parn
\vbox{
$$ \mbox{\textsl{Find}}~
u \in C([0, T), \Ff) \quad \mbox{\textsl{such that}}$$
\beq u(t) = \UU{t}
u_0 + \int_{0}^t~ d s~ \UU{(t - s)} \PPP(u(s),u(s)) \quad
\mbox{\textsl{for all} $t \in [0, T)$}~. \label{int} \feq
}
\noindent
Of course, we say that a solution $u$ is maximal if it cannot
be extended and global if its domain is $[0,+\infty)$.
The following facts are well known. \parn
(i) With a bit more regularity ($u_0 \in \Fp$ and $\Fm$ reflexive), \rref{int} is equivalent
to the Cauchy problem $d u/ d t = \AA u + \PPP(u,u)$, $u(0) = u_0$ for an unknown function
$u \in C([0,T), \Fp) \cap C^1([0,T), \Fm)$. Independently of these stronger assumptions,
a solution $u$ of \rref{int} is usually called a ''mild solution'' of the Cauchy problem. \parn
(ii) \rref{int} has a unique maximal solution (and any other solution is
a restriction of the maximal one); \parn
(iii) \rref{int} has a global solution if the initial datum is small. Among the many references
available on this point, for convenience we
refer to our work \cite{due}, Proposition 5.12, yielding the following statement:
\begin{prop}
\label{xvx}
\textbf{Proposition.} Suppose
\beq 4 K N \| u_0 \| \leqs 1~; \label{coquaaz} \feq
then, the problem \rref{int} has a global solution $u : [0, +\infty) \vain \Ff$. Furthermore,
\beq \| u(t) \| \leqs  \XX(4 K N \| u_0 \|) \, e^{-B t}\, \| u_0 \| \qquad
\mbox{for $t \in [0,+\infty)$}~, \label{xbbz} \feq
where $\XX \in C([0,1], [1,2])$ is the increasing function defined by
\beq \XX(z) := \left\{ \barray{ll} \dd{1 - \sqrt{1 - z} \over (z/2)} & \mbox{for $z \in (0,1]$}~, \\
1 & \mbox{for $z = 0$}~. \farray \right.
\label{xquaazm} \feq
More roughly, due to the features of the function $\XX$, one has
\beq \| u(t) \| \leqs 2 \, e^{-B t}\, \| u_0 \| \qquad
\mbox{for $t \in [0,+\infty)$}~. \feq
\end{prop}
\vskip 0.1cm \noindent
Let us mention that the global solution $u$ in Proposition \ref{xvx} is obtained as the limit
of a Picard iteration: $u(t) = \lim_{k \vain +\infty} u_k(t)$
for each $t \in [0,+\infty)$, where $(u_k)_{k=0,1,2,...}$ is the
sequence of functions in $C([0,+\infty), \Ff))$ defined by $u_0(t) := 0$
and $u_{k+1}(t) := e^{t \AA} u_0 + \int_{0}^t d s \, e^{(t-s) \AA} \PPP(u_k(s),u_k(s))$.
\section{Sobolev spaces and the Na\-vier-Stokes equations on a torus}
\label{inns}
In this section we consider any space dimension
\beq d \geqs 2~. \feq
We use
$r,s$ as indices running from $1$ to $d$;
elements $a, b$,.. of $\reali^d$ or $\complessi^d$ are written with upper or
lower indices, according to convenience: $(a^r)$ or
$(a_r)$, $(b^r)$ or $(b_r)$. For $a$, $b$ $\in \complessi^d$ (say, with upper indices), we put
\beq a \, \sc \, b := \sum_{r=1}^d a^r \, b^r~; \qquad |a| := \sqrt{\overline{a} \, \sc \, a} \feq
where $\overline{a} := (\overline{a^r})$ is the complex conjugate of $a$.
Hereafter we refer to the $d$-dimensional torus
\beq \Td := \underbrace{\To \times ... \times \To}_{\tiny{\mbox{$d$ times}}}~, \qquad \To := \reali/(2 \pi \interi)~, \feq
whose elements are typically written
$x = (x^r)_{r=1,...d}$.
\salto
\textbf{Distributions on $\boma{\Td}$, Fourier series and
Sobolev spaces.}
We introduce the space of periodic distributions $\DD'(\Td, \complessi) \equiv \DD'_{\comple}$, which is
the dual of $C^{\infty}(\Td, \complessi) \equiv C^{\infty}_{\comple}$ (equipping the latter
with the topology of uniform convergence of all derivatives);
$\la v, f \ra$ denotes the action of a distribution $v \in \DD'_{\comple}$
on a test function $f \in C^{\infty}_{\comple}$. \parn
We also consider the lattice $\Zd$ of elements $k =
(k_r)_{r=1,...,d}$ and the Fourier basis $(e_k)_{k \in \Zd}$, where
\beq e_k: \Td \vain \complessi~, \qquad e_k(x) := {1 \over (2 \pi)^{d/2}} \, e^{i k \sc \, x}~
\feq
($k \sc \, x = \sum_{r=1}^d k_r x^r$ makes sense as an element of $\Td$).
Each $v \in \DD'_{\comple}$
has a unique (weakly convergent) Fourier series expansion
\beq v = \sum_{k \in \Zd} v_k e_k~, \qquad
v_k := \la v, e_{-k} \ra~ \in \complessi~. \label{fs} \feq
(As well known, the Fourier coefficients $v_k$ of any $v \in \DD'_{\comple}$
grow polynomially in $|k|$; conversely, any family $(v_k)_{k \in \Z^d}$
of complex numbers with such a polynomial growth is the family of the
Fourier coefficients of some $v \in \DD'_{\comple}$).
The \textsl{mean} of $v \in \DD'_{\comple}$ is
\beq \la v \ra := {1 \over (2 \pi)^d} \la v, 1 \ra = {1 \over (2 \pi)^{d/2}} v_0 \label{mean} \feq
(of course,
$\la v, 1 \ra = \int_{\Td} dx\, v(x)$ if $v \in L^1(\Td)$); the space
of \textsl{zero mean} distributions is
\beq \DD'_{\comple \zz} := \{ v \in \DD'_{\comple}~|~\la v \ra = 0 \}~. \feq
The distributions $v$ in this subspace are characterized by the equivalent condition $v_0 = 0$;
so, their relevant Fourier coefficients are labeled by the set
\beq \Zd_\zz := \Zd \setminus \{0\}~. \feq
The complex conjugate of a distribution $v \in \DD'_{\comple}$ is the unique distribution $\overline{v}$ such that
$\overline{\la v, f \ra} = \la \overline{v}, \overline{f} \ra$ for each $f
\in C^{\infty}_{\comple}$; one has $\overline{v}= \sum_{k \in \Zd} \overline{v_{k}} \, e_{-k}$. \parn
The distributional derivatives $\partial/\partial x^s \equiv \partial_{s}$ $(s=1,...,d)$ and the Laplacian
$\Delta := \sum_{s=1}^d \partial_{s s}$ obviously send $\DD'_{\comple}$ in
$\DD'_{\comple \zz}$, and are such that, for any $v \in \DD_{\comple}$,
$\partial_s v = i \sum_{k \in \Zd_\zz} k_s v_k e_k$,
$\Delta v = - \sum_{k \in \Zd_\zz} | k |^2 v_k e_k$.
For any $n \in \reali$, we further define
\beq \sqrt{-\Delta}^{\, n} : \DD'_{\comple} \vain \DD'_{\comple \zz}~, \qquad
v \mapsto \sqrt{- \Delta}^{\, n}  v := \sum_{k \in \Zd_\zz} | k |^{n} v_k e_k~. \label{reg} \feq
In the sequel, we are interested in the space of
\textsl{real} distributions $\DD'(\Td, \reali) \equiv \DD'$, defined as follows:
\beq \DD' := \{ v \in \DD'_{\comple}~|~\overline{v} = v \} =
\{ v \in \DD'_{\comple}~|~\overline{v_k} = v_{-k}~\mbox{for all $k \in \Zd$} \}~; \feq
of course, $v \in \DD'$ implies $\la v \ra \in \reali$. We also set
\beq \DD'_\zz := \{ v \in \DD'~|~~\la v \ra = 0 \}~; \feq
all the differential operators mentioned before send $\DD'$ into $\DD'_{\zz}$. \parn
Hereafter we consider the real Hilbert space $L^2(\Td, \reali, d x) \equiv L^2$,
with the inner product
$\la v | w \ra_{L^2} := \int_{\Td} v(x) w(x) d x = \sum_{k \in \Zd} \overline{v_k} w_k$ and
the induced norm
\beq \| v \|_{L^2} = \sqrt{\int_{\Td} v^2(x) d x} = \sqrt{\sum_{k \in \Zd} | v_k |^2}~. \label{norms} \feq
To go on, we introduce the \textsl{zero mean Sobolev spaces} $H^n_\zz(\Td, \reali) \equiv H^n_\zz$. For each $n \in \reali$,
\beq H^n_\zz := \{ v \in \DD'_\zz~|~\sqrt{-\Delta}^{\, n} v \in L^2 \} =
\{ v \in \DD'_\zz~|~\sum_{k \in \Zd_\zz} | k |^{2 n} | v_k |^2 < + \infty~\}~;
\label{hn} \feq
this is a real Hilbert space with inner product
$\la v | w \ra_{n} := \la \sqrt{-\Delta}^{\, n} v \, | \, \sqrt{- \Delta}^{\, n} w \ra_{L^2}$
$= \sum_{k \in \Zd_\zz} | k |^{2 n} \, \overline{v_k}  w_k$ and the induced norm
\beq \| v \|_{n} = \| \sqrt{-\Delta}^{\, n}~ v \|_{L^2} = \sqrt{\sum_{k \in \Zd_\zz} | k |^{2 n} | v_k |^2}~.
\label{refur} \feq
Let us consider, in particular, the case $n$ is a \emph{nonnegative integer}. For all $v \in \DD'$ and $k \in \Zd$,
one has $|k|^{2 n} |v_k|^2 = (\sum_{s=1}^d k_s^2)^n |v_k|^2$ $= \sum_{s_1,...,s_n=1}^d k^2_{s_1} ... k^2_{s_n} |v_k|^2$
$= \sum_{s_1,...,s_n=1}^d |(\partial_{s_1 ... s_n} v)_k|^2$; so,
\parn
\vbox{
\beq H^n_\zz := \{ v \in \DD'_\zz~|~\partial_{s_1 ... s_n} v \in L^2~~ \forall s_1,...,s_n \in \{1,....,d\} \}~;
\label{hnn} \feq
\beq \| v \|_n = \sqrt{\sum_{s_1,...,s_n=1}^d \| \partial_{s_1 ... s_n} v \|^2_{L^2}}~
\qquad \forall v \in H^n_\zz~. \label{norintero} \feq
} \noindent
All these details about Sobolev spaces may seem pleonastic, but in fact they
are useful to ensure that some norms we use later coincide exactly
with the norms of \cite{Rob}. \parn
For (integer or noninteger) $n' \leqs n$, one has
$H^{n}_\zz \hookrightarrow H^{n'}_\zz$ and $\|~\|_{n'} \leqs \|~\|_{n}$ on $H^{n}_\zz$.
In particular, $H^\zz_\zz$ is the subspace of $L^2$ made of zero mean elements.
For any real $n$, $\Delta H^n_\zz = H^{n-2}_\zz$ and $\Delta$ is continuous
between these spaces.
Obviously enough, we could define as well the complex Hilbert spaces $L^2_{\comple}$
and $H^n_{\comple \zz}$; however, these are never needed in the sequel. Other facts
about Sobolev spaces (in particular, the duality between $H^{-n}_\zz$ and $H^{n}_\zz$)
are mentioned when necessary in the sequel.
\salto
\textbf{Spaces of vector valued functions on $\Td$.}
If $\Vi(\Td, \reali) \equiv \Vi$ is any vector space of \textsl{real}
functions or distributions on $\Td$, we write
\beq \bVi(\Td) \equiv \bVi := \{ v = (v^1,...,v^d)~|~v^r \in \Vi \quad \mbox{for all $r$}\}~. \feq
In this way we can define, e.g., the spaces $\bb{\DD}'(\Td) \equiv \bb{\DD}'$, $\bb{L}^2(\Td)
\equiv \bb{L}^2$, $\Hz{n}(\Td) \equiv \Hz{n}$.
Any $v = (v^r) \in \bb{\DD}'$ is referred to
as a (distributional) \textsl{vector field} on $\Td$. We note that $v$ has a unique Fourier
series expansion \rref{fs} with coefficients
\beq v_k = (v^r_k)_{r =1,...,d} \in \complessi^d~, \qquad v^r_k  := \la v^r, e_{-k} \ra~; \feq
as in the scalar case, the reality of $v$ ensures $\overline{v_k} = v_{-k}$.
We define componentwisely the mean $\la v \ra \in \reali^d$
of any $v \in \bb{\DD}'$ (see Eq. \rref{mean});
$\Dz$ is the space of zero mean vector fields. We similarly define
componentwisely the operators
$\partial_s, \Delta, \sqrt{-\Delta}^{\,n} : \bb{\DD}' \vain \Dz$.
$\bb{L}^2$ is a real Hilbert space; its
inner product is as in the line before \rref{norms}, with $v(x) w(x)$ and $\overline{v_k} w_k$
replaced by $v(x) \sc \, w(x) = \sum_{r=1}^d v^r(x) \, w^r(x)$
and $\overline{v_k} \sc \, w_k = \sum_{r=1}^d \overline{v^r_{k}} w^{r}_{k}$.
For any real $n$, the $n$-th Sobolev space of zero mean vector fields $\Hz{n}(\Td) \equiv \Hz{n}$ is made
of all $d$-uples $v$ with components
$v^r \in H^n_\zz$; an equivalent definition can be given via Eq.\rref{hn},
replacing therein $L^2$ with $\bb{L}^2$.
$\Hz{n}$ is a real Hilbert space with the inner product
$\la v | w \ra_{n} := \la \sqrt{-\Delta}^{\, n} \, v \, | \sqrt{-\Delta}^{\,n} \,w \ra_{L^2}
= \sum_{k \in \Zd_\zz} | k |^{2 n} \, \overline{v_k}  \sc \, w_k$;
the induced norm $\|~\|_{n}$ is given, verbatim, by Eq. \rref{refur}. If
$w \in \Hz{n}$ has components $w^r$, we have obviously
\beq \| w \|_n = \sqrt{\sum_{r=1}^d \| w^r \|^2_n}~; \feq
for $n$ a nonnegative integer, the vector analogue of Eq. \rref{norintero} is
\beq \| w \|_n = \sqrt{\sum_{r,s_1,...,s_n=1}^d \| \partial_{s_1 ... s_n} w^r \|^2_{L^2}}~.
\label{norvintero} \feq
\salto
\textbf{Divergence free vector fields.} Let us consider the divergence operator
$\dive : \bb{\DD}' \vain \DD'_\zz$, $v \mapsto \dive \,v := \sum_{r=1}^d \partial_r v^r$;
of course,
$\dive \, v = i \, \sum_{k \in \Zd_\zz}  (k \sc \, v_k) e_k$. The
space of \textsl{divergence free (or solenoidal) vector fields} is
\beq \Ds := \{ v \in \bb{\DD}'~|~\dive \,v = 0 \}
= \{ v \in \bb{\DD}'~|~k \sc \, v_k = 0~\forall k \in \Zd~\} ~; \feq
in the sequel, we consider as well the spaces
\beq \Dsz := \Ds \cap \Dz~, \qquad
\HM{n} := \Ds \cap \Hz{n} ~~~(n \in \reali) \feq
(as usually, reference to $\Td$ is omitted for brevity: for example,
$\HM{n}$ stands for $\HM{n}(\Td)$).
$\HM{n}$ is a closed subspace of the Hilbert space $\Hz{n}$,
equipped with the restriction of the inner product $\la~|~\ra_n$. We note the
relations $\Delta(\Ds) = \Dsz$ and $\Delta(\Hs{n+2}) = \HM{n}$, following
immediately from the Fourier representations.
\parn
The \emph{Leray projection} is the map
\beq \LP : \bb{\DD}' \vain \Ds~, \qquad
v \mapsto \LP v := \sum_{k \in \Zd} (\LP_k v_k) e_k~, \feq
where, for each $k$,
$\LP_k$ is the orthogonal projection of $\complessi^d$ onto the orthogonal
complement of $k$;
more explicitly, if $c \in \complessi^d$,
\beq \LP_\zz c = c~, \qquad \LP_k c = c - {k \sc \, c  \over | k |^2}\, k \quad \mbox{for $k \in \Zd_\zz$}~. \feq
From the Fourier representations of $\LP$, $\la~\ra$, $\partial_s$, etc., one easily infers the following statements:
\beq \la \LP v \ra = \la v \ra,\quad
\LP (\partial_s v) = \partial_s (\LP v),\quad \LP (\Delta v) = \Delta (\LP v)\qquad
\forall v \in \bb{\DD}'; \feq
\beq \LP \Dz = \bb{D}'_{\so}, \quad \LP \Hz{n} = \HM{n}; \quad
\| \LP v \|_n \leqs \| v \|_n~ \forall n \in \reali, v \in \Hz{n}~. \label{tvv} \feq
\salto
\textbf{The spaces $\boma{\Hz{1}}$, $\boma{\HM{1}}$.} In the sequel, we are often interested in
specializing the previous considerations to the case $n=1$. This carries
us to the space
$$ \Hz{1} := \{ v \in \Dz~|~\sqrt{-\Delta} \, v \in L^2 \} =
\{ v \in \Dz~|~\partial_s v^r \in L^2 ~\forall r,s \in \{1,...,d\} \} $$
\beq = \{ v \in \Dz~|~\sum_{k \in \Zd_\zz} | k |^{2} | v_k |^2 < + \infty~\}~,
\label{hn1} \feq
with the norm
\beq \| v \|_{1} = \| \sqrt{-\Delta}~ v \|_{L^2} = \sqrt{ \sum_{r,s=1}^d \| \partial_s v^r \|_{L^2}^2} =
\sqrt{\sum_{k \in \Zd_\zz} | k |^{2} | v_k |^2}~.
\label{refur1} \feq
We are as well interested in the divergence free part $\HM{1} = \Ds \cap \Hz{1}$. We note that
\beq \| v \|_1 = \| \curl v \|_{L^2} \qquad \mbox{if $d=3$ and $v \in \HM{1}$}~, \label{cur} \feq
where the right hand side contains the curl of vector fields
(i.e., $(\curl v)^1 := \partial_2 v^3 - \partial_3 v^2$, etc.,
intending derivatives in the distributional sense).
To check Eq. \rref{cur}, one notes that
$(\curl v)_k = i k \wedge v_k$ where $\wedge$ is the (complexified) vector product
(so that $(k \wedge a)^1 := k_2 a^3 - k_3 a^2$, etc., for all $a \in \complessi^3$); if $v$ is
divergence free, the orthogonality property $k \sc v_k = 0$ implies $|k \wedge v_k | = |k | |v_k|$, whence
$\| \curl v \|_{L^2}^2 = \sum_{k \in \Z^3_\zz} |k|^2 |v_k|^2 = \| v \|^2_1$. \parn
(Before going on, we point out that the norm
$\| v \|_1$ on $\HM{1}$ (in dimension $d=3$) is denoted in
\cite{Rob} with $\| D v \|$, the symbol $\|~\|$ standing for the $L^2$ norm:
see page 40 of the cited reference. The equality $\| v \|_1 = \| D v \|$ is essential
for subsequent comparison between some estimates of ours and \cite{Rob}).
\salto
\textbf{The exponential of the Laplacian}. Let us put
\beq e^{t \Delta} v := \sum_{k \in \Zd} e^{-t | k |^2} v_k e_k~ \qquad \forall t \in [0,+\infty),
v \in \bb{D}'~; \label{etd} \feq
the above series converges in the weak topology of $\bb{D}'$. \parn
From the Fourier representations, it is clear that
\beq e^{t \Delta} \Dz \subset \Dz~, \qquad e^{t \Delta} \Dsz \subset \Dsz
\quad \forall t \in [0,+\infty). \feq
In the sequel, we consider any $n \in \reali$.
As well known, the map $(t, v) \mapsto e^{t \Delta} v$
sends continuously $[0,+\infty) \times \Hz{n}$ into $\Hz{n}$; this map
is a strongly continuous semigroup on the Hilbert space $\Hz{n}$, with generator
$\Delta \restriction \Hz{n+2} : \Hz{n+2} \vain \Hz{n}$. Furthermore
\beq  \no{e^{t \Delta} v}_n \leqs e^{-t} \no{v}_n
\qquad \forall t \in [0,+\infty),
v \in \Hz{n} \label{eqn} \feq
(as follows immediately from the Fourier representations
\rref{etd} of $e^{t \Delta}$, combined with the obvious inequality
$e^{- t |k |^2} \leqs e^{-t}$ for all $k \in \Zd_\zz$). \parn
Now, let us consider two Sobolev spaces $\Hz{n}$, $\Hz{n-\nu}$, with $\nu \in (0,+\infty)$; then
\beq e^{t \Delta} v \in \Hz{n},
\qquad \no{e^{t \Delta} v}_n \leqs \mmu_{\nu}(t) e^{-t} \no{v}_{n - \nu}
\qquad \forall t \in (0,+\infty), v \in \Hz{n- \nu}, \label{eqvun} \feq
\beq
\mmu_\nu(t) := \left\{ \barray{ll} \dd{\left({\nu \over {2 e t}}\right)^{{\nu \over 2}} e^t}
& \mbox{for $0 < t \leqs \dd{{\nu \over 2}}$}~,
\\ 1 & \mbox{for $t > \dd{ {\nu \over 2}}$}~; \farray \right. \label{eqwun} \feq
(see Appendix \ref{appesp}). The map $(t, v) \mapsto e^{t \Delta} v$
sends continuously $(0,+\infty) \times \Hz{n - \nu}$ into $\Hz{n}$.
\parn
All the previous statements about $\dd{e^{t \Delta}}$ and the Sobolev spaces, including
the relations \rref{eqn}-\rref{eqwun}, hold
as well if we replace systematically $\Hz{n}$, $\Hz{n+2}$, $\Hz{n-\nu}$ with $\HM{n}$, $\HM{n+2}$,
$\HM{n-\nu}$; this is basically due to the relations $\Delta \Dsz = \Dsz$ and $e^{t \Delta}
\Dsz \subset \Dsz$. \salto
\textbf{NS equations.} A setting for the incompressible NS equations on $\Td$
consists of three Banach spaces $\Fp, \Ff, \Fm$ of divergence free vector fields on $\Td$,
such that all conditions (\PU-\Q) of Section \ref{abf} are fulfilled by these spaces, by the operator
\beq \AA := \Delta \restriction \Fp : \Fp \mapsto \Fm \label{defa}\feq
and by the bilinear map
\beq \PPP : \Ff \times \Ff \vain \Fm~, \qquad (v, w) \mapsto \PPP(v,w) := - \LP(v \sc \partial w)~. \feq
In the above, $v \sc \partial w$ is the vector field on $\Td$ with components $(v \sc \partial w)^r :=
\sum_{s = 1}^d v^s \partial_s w^r$; the function spaces are chosen so that,
for $v, w \in \Ff$, the pointwise products $v^s \partial_s w^r$ are well defined  and $\PPP(v, w)$
belongs to $\Fm$. \parn
The initial value problem \rref{int} corresponding to the above maps $\AA, \PPP$ and to any datum $u_0 \in \Ff$
takes the form
$$ \mbox{\textsl{Find}}~
u \in C([0, T), \Ff) \quad \mbox{\textsl{such that}}$$
\beq u(t) = e^{t \Delta}
u_0 - \int_{0}^t~ d s~ e^{(t-s) \Delta} \LP(u(s) \sc \partial u(s)) \quad
\mbox{\textsl{for all} $t \in [0, T)$}~, \label{intns} \feq
and is related, in the sense already indicated, to the Cauchy problem
$d u/ d t = \Delta u - \LP(u \sc \partial u)$, $u(0) = u_0$. One recognizes the
NS evolution equation for an incompressible fluid (in units in which the
density and the viscosity are $1$, and assuming no external forces). \parn
In \cite{due} we considered for arbitrary $d$ the setting
$\Ff := \HM{n}$, $\Ff_{\pm} := \HM{n\pm 1}$ $(n \in (d/2, + \infty))$,
putting the emphasis on all the quantitative aspects and, in particular,
on the accurate estimation of the constants $N, K$ in Section \ref{abf}; this allowed
to get fully quantitative estimates on the time of existence for the
solution of \rref{int}, and gave a framework to evaluate \textsl{a posteriori} the distance between any approximate
solution and the exact solution. \parn
As anticipated in the Introduction, here we will take a similar attitude for the setting
\beq \Fp := \HM{2 - \del}~, \quad \Ff := \HM{1}~, \quad \Fm := \HM{-\del}~,  \label{define} \feq
mainly for $d=3$ and $1/2 < \del < 1$.  Again, our evaluation of $N,K$ will be
relevant in connection with the time of existence for \rref{intns}, with the
estimates about approximate solutions, and so on.
\section{An $\boma{H^1}$ setting for the semigroup of the Laplacian}
\label{setting}
For the moment we assume
\beq d \in \{2,3,....\}~, \qquad \del \in (0,1) \feq
and take  $\Fm$, $\Ff$, $\Fp$ as in \rref{define},
$\AA$ as in \rref{defa}; the inclusion (\PU) $\Fm \hookrightarrow \Ff \hookrightarrow \Fp$ is evident.
\begin{prop}
\label{edelta}
\textbf{Proposition.}
(i) $\AA$ generates a strongly continuous semigroup on $\HM{-\del} = \Fm$; this is
the restriction to $[0,+\infty) \times \HM{-\del}$ of the map $(t,v) \mapsto e^{t \Delta} v$
(see Eq. \rref{etd}). \parn
(ii) The map $(t,v) \mapsto e^{t \Delta} v$ gives as well a strongly continuous semigroup on
$\HM{1} = \Ff$; one has
\beq \no{e^{t \Delta} v}_1 \leqs e^{-t} \no{v}_1
\qquad \forall t \in [0,+\infty), v \in \HM{1}~. \label{reguff} \feq
(iii) The map $(t, v) \mapsto e^{t \Delta} v$ sends continuously
$(0,+\infty) \times \HM{-\del} = (0,+\infty) \times \Fm$ into $\HM{1} = \Ff$;
one has
\beq \no{e^{t \Delta} v}_1 \leqs \mu_{\del}(t) e^{-t} \no{v}_{-\del}
\qquad \forall t \in (0,+\infty), v \in \HM{-\del}~, \label{eqvv} \feq
\beq
\mu_\del(t) := \left\{ \barray{ll} \dd{\left({1 + \del \over {2 e t}}\right)^{{1 + \del \over 2}} e^t}
& ~~\mbox{for~~ $0 < t \leqs \dd{{1 + \del \over 2}}$}~,
\\ 1 & ~~\mbox{for ~~$t > \dd{ {1 + \del \over 2}}$}~. \farray \right. \label{eqww} \feq
(Note that $\mu_\del(t) = O(1/t^{(1 + \del)/2})$ for $t \vain 0^+$; this power
of $1/t$ is integrable due to the assumption $0 < \del < 1$).
With $\mu_\del$ as above, one has
\beq \sup_{t \in [0,+\infty)}
\int_{0}^t d s \,e^{-s}\, \mu_\del(t-s)  < + \infty~. \label{defn} \feq
(iv) Statements (i-iii) indicate that
the set $(\Fp, \Ff, \Fm, \AA)$ fulfills all conditions (\PD-\PC) of
Section \ref{abf}. The function $\mu$ and the constants $B, N \in (0,+\infty)$
mentioned therein can be taken as follows:
\beq \mu = \mu_{\del}~\mbox{as in \rref{eqww}}, \quad B = 1, \quad N \equiv N_{\del} := \mbox{any majorant
of the sup \rref{defn}}. \label{holdw} \feq
\end{prop}
\textbf{Proof.} (i-iii) Essentially, one specializes to the present case
all statements of Section \ref{inns} (and Appendix \ref{appesp})
about $(e^{t \Delta})$ and Sobolev spaces of arbitrary order.
In particular, the relations \rref{reguff}-\rref{eqww} follow
from the relations \rref{eqn}-\rref{eqwun} for the spaces
$\HM{n}$ and $\HM{n-\nu}$, with $n=1$ and $\nu = 1 + \del$. With
this choice of $\nu$, the function
$\mmu_{\nu}$ of Eq. \rref{eqwun} becomes the $\mu_{\del}$ of
Eq. \rref{eqww}; the finiteness of the sup \rref{defn} is easily inferred from
the explicit expression of $\mu_{\del}$, by an argument very similar to the one
employed in \cite{uno} to prove Proposition 7.2, item (iv). \parn
(iv) Obvious. \fine
\section{An $\boma{H^1}$ setting for the NS quadratic nonlinearity}
\label{sequad}
Throughout this section we assume
\beq d \in \{2,3,...\}~, \qquad \del \in ({d \over 2} -1, +\infty)~. \label{assd} \feq
Our aim is to show that the NS bilinear map
$(v,w) \mapsto \PPP(v, w) :=-  \LP(v \sc \partial w)$
is well defined from $\Ff \times \Ff  = \HM{1} \times \HM{1} $ to $\Fm = \HM{-\del}$,
and to determine a constant
$K_{\del}$ such that $\| \PPP(v, w) \|_{-\del} \leqs K_{\del} \| v \|_1 \| w \|_1$. \parn
We will reach this result through a number of intermediate steps, most of them
relying on the function
\beq \KK_{\del} : \Zd_\zz \vain (0,+\infty)~, \qquad k \mapsto \KK_{\del}(k) :=
 \sum_{h \in \Zd_{0 k}} {1 \over |h|^{2 \del} \, |k - h|^2}~;
\label{kapde} \feq
in the above one should intend
\beq \Zd_{0 k} := \Zd \setminus \{0,k\}~, \feq
a notation to be employed systematically in the sequel. \parn
In Appendix \ref{appeka} we show the following.\parn
(i) The sum defining $\KK_{\del}(k)$ (in principle, existing
in $(0,+\infty]$) is in fact finite, for each $k \in \Zd_\zz$. \parn
(ii) One has $\sup_{k \in \Zd_\zz} \KK_{\del}(k) < + \infty$. \parn
(iii) Explicit lower and upper bounds can be given for $\KK_{\del}$, to be used
for the practical evaluation of this function and of its $\sup$.
\vskip 0.2cm \noindent
The function $\KK_{\del}$ first appears in the forthcoming Proposition,
dealing with pointwise multiplication in the Sobolev spaces
(of \emph{scalar} functions) $H^{\del}_\zz$, $H^1_\zz$; this Proposition will be subsequently employed to estimate the
NS bilinear map.
\begin{prop}
\label{kadel}
\textbf{Proposition.} With $d$ and $\del$ as in \rref{assd}, the following holds. \parn
(i) Let $\z \in H^{\del}_\zz$, $v \in H^{1}_\zz$, so that $\z v$ (the pointwise product of two
$L^2$ functions) is well defined (and in $L^1$). This product has the additional features
\beq \z v \in L^2~, \qquad \| \z v - \la \z v \ra \|_{L^2} \leqs K_{\del} \| \z \|_{\del} \| v \|_{1}  \label{kuv} \feq
for a suitable constant $K_{\del} \in (0,+\infty)$, independent of $z, v$ (recall that $\la~\ra$
is the mean). \parn
(ii) The above inequality is fulfilled for all $\z, v$ by any constant $K_{\del}$ such that
\beq {1 \over (2 \pi)^{d/2}} \sqrt{\sup_{k \in \Zd_{0}} \KK_{\del}(k)} \leqs K_{\del}, \label{dekde} \feq
with $\KK_{\del}$ as in Eq. \rref{kapde}.
\end{prop}
\textbf{Proof.} Let us consider the Fourier coefficients $(\z v)_k$ of $\z v$, for
$k \in \Zd$;  up to a factor $1/(2 \pi)^{d/2}$, these are obtained taking the
convolution of the Fourier coefficients of $\z,v$:
\beq (\z v)_k = {1 \over (2 \pi)^{d/2}} \sum_{h \in \Zd} \z_{h} v_{k - h}~ =
{1 \over (2 \pi)^{d/2}} \sum_{h \in \Zd_{0 k}} \z_{h} v_{k -h}~  \feq
(the sum can be confined to $\Zd_{0 k}$, due to the vanishing of $\z_0, v_0$).
From here we get
\beq |(\z v)_k| \leqs {1 \over (2 \pi)^{d/2}} \sum_{h \in \Zd_{0 k}} |\z_{h}| | v_{k - h} |  \feq
$$ = {1 \over (2 \pi)^{d/2}} \sum_{h \in \Zd_{0 k}}~{1 \over |h|^\del |k - h|}~
|h|^\del | \z_{h} | | k - h |~| v_{k - h} |~. $$
Now, H\"older's inequality $| \sum_h~ a_h b_h |^2 \leqs \Big(\sum_h | a_h |^2\Big)
\Big(\sum_h~| b_h |^2 \Big)$ gives
\beq | (\z v)_k |^2 \leqs {1 \over (2 \pi)^d} \KK_\del(k) \PP(k)~, \label{dains} \feq
$$ \KK_\del(k) := \sum_{h \in \Zd_{0 k}} {1 \over |h|^{2 \del} |k - h|^2}~
\mbox{as in \rref{kapde}}, \quad
\PP(k) := \sum_{h \in \Zd_{0 k}} |h|^{2 \del} |  \z_{h} |^2
| k - h |^2 | v_{k - h} |^2~. $$
Eq. \rref{dains} implies
$$ \sum_{k \in \Zd_{0}} |(\z v )_k|^2 \leqs
{1 \over (2 \pi)^d} \sum_{k \in \Zd_{0}} \KK_\del(k) \PP(k) $$
\beq \leqs {1 \over (2 \pi)^d} \Big( \sup_{k \in \Zd_{0}} \KK_\del(k) \Big) \,
\sum_{k \in \Zd_0} \PP(k) \leqs K_{\del}^2 \sum_{k \in \Zd_0} \PP(k)
\label{dains0} \feq
where $K_{\del}$ is any constant such that $(2 \pi)^{-d/2}
\sqrt{\sup_{k \in \Zd_{0}} \KK_\del(k)}$ $\leqs K_{\del}$, as in \rref{dekde}. \parn
To continue, let us observe that the definition of $\PP$ implies
$$ \sum_{k \in \Zd_{0}} \PP(k) = \sum_{(k, h) \in \Zd_0 \times \Zd_0, k \neq h} |h|^{2 \del} |  \z_{h} |^2
| k - h |^2 | v_{k - h} |^2 $$
whence, with a change of variable $\ell = k - h$,
\beq \sum_{k \in \Zd_{0}} \PP(k) = \sum_{(h, \ell) \in \Zd_0 \times \Zd_0} |h|^{2 \del} |  \z_{h} |^2
| \ell |^2 | v_{\ell} |^2 =
\| \z \|^2_\del \, \| v \|^2_1~.
\label{dains2} \feq
From Eqs.\rref{dains0} \rref{dains2}, we get
\beq \sum_{k \in \Zd_\zz} |(\z v )_k|^2 \leqs K_{\del}^2 \| \z \|^2_\del \, \| v \|^2_1~. \label{laseq} \feq
The finiteness of $\sum_{k \in \Zd_\zz} |(\z v )_k|^2$, and thus of $\sum_{k \in \Zd} |(\z v )_k|^2$,
ensures $\z v \in L^2$. Noting that $\z v - \la \z v \ra = \sum_{k \in \Zd_\zz} (\z v )_k e_k$ we can
reexpress \rref{laseq} as
\beq \| \z v - \la \z  v \ra \|^2_{L^2} \leqs K_{\del}^2 \| \z \|^2_\del \, \| v \|^2_1~, \feq
whence the thesis. \fine
\begin{prop}
\label{ppp}
\textbf{Proposition.} Let $d$, $\del$ be as in \rref{assd}, and denote with $K_{\del}$
any constant fulfilling Eq. \rref{dekde}. Furthermore, let
\beq v \in \HM{1}~, \qquad w \in \Hz{1}~; \feq
then (i)(ii) hold. \parn
(i) Consider the product $v \sc \partial w$ (whose components
$(v \sc \partial w)^r :=\sum_{s=1}^d v^s \partial_s w^r$ are well defined and belong to $L^1$,
being sum of products of the $L^2$ functions $v^s$, $\partial_s w^r$). One has
\beq v \sc \partial w \in \Hz{-\del}~, \qquad \| v \sc \partial w \|_{-\del} \leqs K_{\del} \| v \|_1
\| w \|_1 ~. \label{evw} \feq
(ii) Furthermore, one has
\beq \LP(v \sc \partial w) \in \HM{-\del}~, \qquad \| \LP(v \sc \partial w) \|_{-\del} \leqs K_{\del} \| v \|_1
\| w \|_1~.\feq
(Due to the above results, condition (\Q) of Section \ref{abf} is fulfilled
with $\Ff := \HM{1}$, $\Fm := \HM{-\del}$, $\PPP : \Ff \times \Ff \vain \Fm$,
$\PPP(v, w) := -\LP(v \sc \partial w)$ and $K := K_{\del}$).
\end{prop}
\textbf{Proof.} Let us observe that (ii) follows immediately from (i) using \rref{tvv} with $n \vain -\del$
and $v \vain v \sc \partial w$. Hereafter, we give the proof of (i); after putting for brevity
\beq \p := v \sc \partial w~, \feq
we proceed in several steps. \parn
\textsl{Step 1. One has
\beq \p^r \in L^1_\zz := L^1 \cap \DD'_\zz \feq
for all $r \in \{1,...,d\}$ (i.e., $\p \in \mathbb{L}^1_\zz := \mathbb{L}^1 \cap \Dz$)}.
We have just noticed, as a comment in the statement of this Proposition, that $\p^r \in L^1$;
furthermore, $\p^r$ has zero mean since
\beq \int_{\Td} \p^r d x = \sum_{s=1}^d \int_{\Td} v^s \partial_s w^r = -
\sum_{s=1}^d \int_{\Td} (\partial_s v^s) w^r = 0 \feq
(here, we have used an integration by parts and the assumption $\dive v = 0$). \parn
\textsl{Step 2. For all $r \in \{1,...,d\}$ one has}
\beq \p^r \in H^{-\del}_\zz~, \qquad
\| \p^r \|_{-\del} \leqs K_{\del} \| v \|_1 \sqrt{\sum_{s=1}^d \| \partial_s w^r \|^2_{L^2}}~.
\label{boundw} \feq
We know that $\p^r \in L^1_\zz$; by the familiar duality between $H^{\del}_\zz$ and $H^{-\del}_\zz$,
\beq \p^r \in H^{-\del}_\zz \qquad \Longleftrightarrow \qquad  \z \p^r \in L^1~ \forall \z \in H^{\del}_\zz,~
\sup_{\z \in H^{\del}_\zz, \z \neq 0} { | \int_{\Td} \z \p^r d x | \over \| \z \|_{\del}} < + \infty~; \label{cr1} \feq
if the above conditions are fulfilled, we further have
\beq \| \p^r \|_{- \del} =  \sup_{\z \in H^{\del}_\zz, \z \neq 0} { | \int_{\Td} \z \p^r d x | \over \| \z \|_{\del}}~.
\label{cr2}\feq
Keeping in mind the above statements, we consider any $\z \in H^{\del}_\zz$ and the function
\beq \z \p^r  = \sum_{s=1}^d \z \, v^s \partial_s w^r~.\feq
For each $s$ we have the following: \parn
(i) $\partial_s w^r \in L^2$, due to the assumption $w \in \Hz{1}$; \parn
(ii)  $\z v^s \in L^2$: this follows using Proposition \ref{kadel}, which also gives
\beq \| \z v^s - \la \z v^s \ra \|_{L^2} \leqs K_{\del} \| \z \|_{\del} \| v^s \|_{1}~. \label{fm0} \feq
From $\z v^s, \partial_s w^r \in L^2$ for all $s$, it follows
\beq \z \p^r \in L^1~; \feq
now we estimate the integral of $\z \p^r$, starting from the following remark:
\beq \int_{\Td} \z \p^r d x = \sum_{s=1}^d \int_{\Td} (\z v^s) \partial_s w^r  \label{fm1} \feq
$$ = \sum_{s=1}^d \int_{\Td} (\z v^s - \la \z v^s \ra) \partial_s w^r d x +
\sum_{s=1}^d \la \z v^s \ra \int_{\Td} \partial_s w^r d x =
\sum_{s=1}^d \int_{\Td} (\z v^s - \la \z v^s \ra) \partial_s w^r~$$
($\int_{\Td} \partial_s w^r d x =0$, since this is the integral of a derivative). From \rref{fm1},
from the H\"older inequality and \rref{fm0} we get
$$  \left| \int_{\Td} \z \p^r d x \right| \leqs \sum_{s=1}^d \| \z v^s - \la \z v^s \ra \|_{L^2} \| \partial_s w^r~\|_{L^2}
\leqs K_{\del} \| \z \|_{\del} \sum_{s=1}^d  \| v^s \|_{1}~\| \partial_s w^r~\|_{L^2} $$
\beq\leqs K_{\del} \| \z \|_{\del}  \sqrt{ \sum_{s=1}^d \| v^s \|^2_{1} }~
\sqrt{ \sum_{s=1}^d \| \partial_s w^r~\|^2_{L^2}}
 = K_{\del} \| \z \|_{\del} \| v \|_{1}~
\sqrt{ \sum_{s=1}^d \| \partial_s w^r~\|^2_{L^2}}~. \label{lam} \feq
Now, using \rref{lam} with \rref{cr1} \rref{cr2} we conclude that $\p^r$ is actually
in $H^{-\del}_\zz$, and $\| \p^r \|_{-\del}$ admits the bound \rref{boundw}.
\parn
\textsl{Step 3. One has}
\beq \p \in \Hz{-\del}~, \qquad \| \p \|_{-\del} \leqs K_{\del} \| v \|_1
\| w \|_1  \label{ezvw} \feq
\textsl{(so, \rref{evw} holds and the proof is concluded)}.
From Step 2, we know that each component $\p^r$ is in $H^{-\del}_\zz$; furthermore,
the estimates \rref{boundw} imply
\beq \| \p\|^2_{- \del} =  \sum_{r=1}^d
\| \p^r \|^2_{-\del} \leqs K^2_{\del} \| v \|^2_1 \sum_{r,s=1}^d \| \partial_s w^r \|^2_{L^2} =
K^2_{\del} \| v \|^2_1 \| w \|^2_{1}~,
\label{boundz} \feq
whence the thesis. \fine
\section{Putting things together: an $\boma{H^1}$ setting for the NS initial value problem,
in dimension $\boma{d=2,3}$}
\label{toge}
Let us recall that
\beq \Fp := \HM{2 - \del}, \quad \Ff := \HM{1}, \quad \Fm := \HM{-\del}~; \qquad
\AA := \Delta \restriction \HM{2-\del} : \HM{2 - \del} \vain \HM{-\del}~; \label{summin} \feq
$$ \PPP : \Ff \times \Ff \vain \Fm ~, \qquad (v, w) \mapsto \PPP(v, w) := - \LP(v \sc \partial w)~. $$
The above Sobolev spaces live, for the moment, on a torus $\Td$ of any dimension $d \in \{2,3,..\}$;
our previous results about this setting can be summarized as follows. \parn
Condition (\PU) of Section \ref{abf} is fulfilled by the triple $\Fp, \Ff, \Fm$;
according to Proposition \ref{edelta}, $\AA$
fulfills conditions (\PD-\PC) if
$\del \in (0,1)$. On the other hand, according to Proposition \ref{ppp},
$\PPP$ is well defined and fulfills condition (\Q) of Section \ref{abf} if
$\del \in (d/2 -1, + \infty)$; so, the requirements of
both Propositions \ref{edelta}, \ref{ppp} hold simultaneously if
\beq {d \over 2} - 1 < \del < 1~. \label{lastc} \feq
Of course, \rref{lastc} can be fulfilled for some $\del$ only if $d/2 - 1 < 1$, i.e., for $d=2$ or $d=3$;
so, \rref{lastc} holds if either
\beq d= 2, \qquad \del \in (0,1) \label{d2} \feq
or
\beq d =3, \qquad \del \in ({1 \over 2}, 1)~. \label{d3} \feq
Summing up (and recalling the statements of Propositions \ref{edelta}, \ref{ppp} on the
function $\mu$ and the constants $B,N,K$), we have the following result.
\begin{prop}
\label{propfin}
\textbf{Proposition.}
Let $d, \del$ be as in Eqs. \rref{d2} \rref{d3}. Then, the set $(\Fp, \Ff, \Fm, \AA, \PPP)$
defined by Eq. \rref{summin} fulfills all conditions (\PU-\Q) of Section \ref{abf}. The function $\mu$
and the constants
$B,N,K$ mentioned in Section \ref{abf} can be taken as follows:
\beq \mu = \mu_{\del}~\mbox{as in \rref{eqww}},~~B = 1,~~N = N_{\del}~\mbox{as in \rref{holdw}},~~
K = K_{\del}~\mbox{as in \rref{dekde}}~. \feq
\end{prop}
The previous proposition allows to apply the full machinery of \cite{due}
to the NS initial value problem \rref{intns}. Hereafter we consider, in particular, the condition for
global existence presented in \cite{due}
and summarized in Proposition \ref{xvx} of the present work.
\salto
\textbf{Global existence for NS with small initial data, when $\boma{d=3}$.}
Let us keep the definitions \rref{summin}.
Global existence for the NS initial value problem \rref{intns} is well known
for any datum in $\Ff = \HM{1}$, if $d=2$; so, we pass to the case $d=3$. Let us
recall that, in the three dimensional case, we have the equality
\rref{cur} $\| v \|_1 = \| \curl v \|_{L^2}$ for all $v \in \HM{1}$.
On the grounds of Propositions \ref{xvx} and \ref{propfin}, we can state the following.
\begin{prop}
\label{xvxd}
\textbf{Proposition.} Let $d =3$, $\del \in ({1 / 2}, 1)$ as in \rref{d3}, and
$\Ff$, etc., as in \rref{summin}. The solution $u$ of the initial value problem \rref{intns} is global
if the initial datum $u_0 \in \Ff = \HM{1}$ is such that
\beq \| u_0 \|_1 \leqs {1 \over 4 N_{\del} K_{\del}}~; \label{coquaazde} \feq
furthermore, the solution fulfills the bound (of the type \rref{xbbz})
\beq \| u(t) \|_1 \leqs  \XX(4 K_\del N_\del \| u_0 \|_1) \, e^{-t}\, \| u_0 \|_1
\qquad \mbox{for $t \in [0,+\infty)$}~, \label{xbbzde} \feq
where $\XX \in C([0,1], [1,2])$ is the increasing function defined by \rref{xquaazm}. More roughly, we have
\beq \| u(t) \|_1 \leqs  2 \, e^{-t}\, \| u_0 \|_1 \qquad \mbox{for $t \in [0,+\infty)$}~. \label{rough} \feq
\end{prop}
To get a fully quantitative estimate, let us put
\beq \del := \delnn \feq
(see the forthcoming Remark \ref{remadel} (iii) about this choice).
Then, computing numerically the function of $t$ indicated below, we see that
\beq \sup_{t \in [0,+\infty)}
\int_{0}^t d s \,e^{-s}\, \mu_{\deln}(t-s)  < 1.70 := N_{\deln} \label{defn34} \feq
(in fact the above sup is a maximum, attained at point $t \in (0.21, 0.22)$; for these
and other numerical computations one can use, e.g., the MATHEMATICA package). \parn
Furthermore, one has (see Appendix \ref{appesup})
\beq 27.94 < \sup_{k \in \Zd_\zz} \KK_{\deln}(k) < 32.23 \label{bounds} \feq
(here we are mainly interested in the upper bound $32.23$;
the lower bound $27.94$ is reported only for an appreciation of
our uncertainty about $\sup \KK_{\deln}$). This implies
\beq 0.335 < {1 \over (2 \pi)^{3/2}} \sqrt{\sup_{k \in \Zd_\zz} \KK_{\deln}(k)} <
0.361 := K_{\deln}~. \feq
The values of
$N_{\deln}$ and $K_{\deln}$ imply
\beq {1 \over 4 N_{\deln} K_{\deln}} > \nume~; \feq
In conclusion, Proposition \ref{xvxd} with the choice $\del=\deln$ and the previous evaluations
of the related constants yield the following.
\begin{prop}
\textbf{Corollary.} With $d=3$,
the solution $u$ of the initial value problem \rref{intns} is global
if the initial datum $u_0 \in \Ff = \HM{1}$ is such that
\beq \| u_0 \|_1 \leqs \nume~; \label{coquaazdet} \feq
furthermore, for all $t \in [0, +\infty)$, the solution fulfills the bounds
\beq \| u(t) \|_1 \leqs  \XX\Big({ \| u_0 \|_1 \over \nume}\Big) \, e^{-t}\, \| u_0 \|_1
\leqs 2 e^{-t } \| u_0 \|_1~,  \label{xbbzdet} \feq
with $\XX$ as in \rref{xquaazm}.
\end{prop}
\begin{rema}
\label{remadel}
\textbf{Remarks.} (i) According to the remark following Proposition \ref{xvx}, the global solution $u$
can be constructed by a Picard iteration. A fine analysis of the iteration, based on the
the smoothing properties of $(e^{t \Delta})$, shows that the function $(t,x) \mapsto u(t)(x) \equiv u(x,t)$
is in fact $C^\infty$ on $(0,+\infty) \times \Tt$, where it satisfies the NS equations in the classical sense;
on this point see, e.g., the proof of Theorem 15.2(A) in \cite{Lem}. \parn
Taking into account these facts and the exponential time decay of $\| u(t) \|_1$, we
see that $u$ is in $L^{\infty}([0,+\infty),  \bb{L}^2_{\so}) \cap L^2([0,+\infty), \HM{1})$
and fulfills the NS equations in the weak sense; so, $u$ is a weak solution of the NS
equations in the sense of \cite{Tem}. Again by the exponential time decay,
it is $u \in L^4([0,+\infty), \HM{1})$; this suffices to infer that $u$ is
a strong NS solution in the sense of \cite{Tem} (see the Remark 3.3 on pages 22-23 of this reference). \parn
(ii) We already mentioned that our criterion \rref{coquaazdet}
improves the condition for the existence of a global (strong) solution
recently proposed in \cite{Rob}; in our notations, the condition of the cited paper reads
\beq \| u_0 \|_1 \leqs \numerob \label{072} \feq
(see page 43 of \cite{Rob}; $\numerob$ is the approximation with three meaningful digits
of the quantity indicated therein by $R_V$). \parn
(iii) Of course, in place of $\del=\deln$ one could consider for $\del$ other choices
in the interval $(1/2,1)$, each one yielding a bound of the type \rref{coquaazde} for
global existence; the best estimate of this type would be
obtained maximizing $1/(4 N_{\del} K_{\del})$, for $\del$ in $(1/2,1)$. However, a few experiments we did with
$\del \neq \deln$
seem to exclude a significant improvement of the bound \rref{coquaazdet}.
\end{rema}
\section{Some possible developments.}
\label{develop}
In \cite{due}, we presented the following abstract result:
\begin{prop}
\textbf{Proposition.}
\label{propapp}
Assume the following: \parn
(i) $(\Fp, \Ff, \Fm, \AA, \PPP)$ is a set with the properties (\PU-\Q)~ of Section
\ref{abf}; $\mu,B,N,K$ are the function and the constants mentioned therein. \parn
(ii) $u_0 \in \Ff$ is a datum for the initial value problem \rref{int}.\parn
(iii) $u_{ap}$ is an approximate solution of problem \rref{int} with domain
$[0,T)$ ($T \in (0,+\infty]$), a growth estimator $\Dd$ and an error estimator $\EE$: these expressions indicate
three functions $u_{ap} \in C([0,T), \Ff)$ and $\Dd, \EE \in C([0,T), [0,+\infty))$ such that
\beq \| u_{ap}(t) \| \leqs \Dd(t)~, \feq
\beq \| u_{ap}(t) - \UU{t} \, u_0 - \int_{0}^t~ d s~ \UU{(t - s)} \PPP(u_{ap}(s),u_{ap}(s))  \| \leqs \EE(t)
\qquad \mbox{for $t \in [0,T)$}~.
\feq
(iv) There is a function $R \in C([0,T), [0,+\infty))$ fulfilling the ''control inequality''
\beq \EE(t) + K \int_{0}^t ds \, \mu(t-s) e^{-B(t-s)} \Big(2 \Dd(s) R(s) + R^2(s)\Big)
\leqs R(t)~\mbox{for $t \in [0,T)$}~. \label{cont}
\feq
Then, problem \rref{int} has an (exact) solution $u \in C([0,T), \Ff)$, and
\beq \| u(t) - u_{ap}(t) \| \leqs R(t) \qquad \mbox{for $t \in [0,T)$}~. \feq
\end{prop}
The previous proposition allows to make predictions on the interval of existence
of the solution $u$ of \rref{int}, and on its distance from $u_{ap}$,
using some information ($\Dd$ and $\EE$)
pertaining to $u_{ap}$ only; in this sense, we have an estimate for $u$ from
an \textsl{a posteriori} analysis of $u_{ap}$. \parn
In \cite{due}, we presented some applications of this result to the NS equations
on $\Td$, taking for $\Ff_{\pm}, \Ff$ some Sobolev spaces of sufficiently high order and
considering, for example, the Galerkin approximate solutions. \parn
The analysis performed in the present work allows to apply Proposition \ref{propapp} to the NS equations with
$\Fp := \HM{2 - \del}, \Ff := \HM{1}, \Fm := \HM{-\del}$ and ($d=2$, $\del \in (0,1)$ or)
$d=3$, $\del \in (1/2,1)$. \parn
Our estimates on some related functions and constants (e.g., $\mu = \mu_{\del}$ and $K = K_{\del}$
considered in Sections \ref{setting}-\ref{toge}) allow a fully quantitative implementation
of the method. It should be pointed out that the availability of accurate information on
$\mu, B, K$, etc., is essential for an efficient application of the control inequality
\rref{cont}; for example, the time interval $[0,T)$ on which a solution $R$ of \rref{cont}
exists depends sensibly on these data. We plan to illustrate elsewhere some applications of
Proposition \ref{propapp} to the NS equations on $\Tt$, in the framework of the above spaces
$\Ff_{\pm}, \Ff$.
\appendix
\section{Appendix. On the semigroup of the Laplacian}
\label{appesp}
After recalling the definition \rref{etd}
$$ e^{t \Delta} v := \sum_{k \in \Zd} e^{-t | k |^2} v_k e_k~ \qquad \forall t \in [0,+\infty),
v \in \bb{D}'~, $$
let us fix $n \in \reali$, $\nu
\in (0,+\infty)$. Hereafter, we prove the following result.
\begin{prop} \textbf{Proposition.} Let $t \in (0,+\infty)$, $v \in \Hz{n-\nu}$. Then,
we have the relations \rref{eqvun}, \rref{eqwun}
$$ e^{t \Delta} v \in \Hz{n}~, \qquad \no{e^{t \Delta} v}_n \leqs \mmu_{\nu}(t) e^{-t} \no{v}_{n - \nu} ~, $$
$$
\mmu_\nu(t) := \left\{ \barray{ll} \dd{\left({\nu \over {2 e t}}\right)^{{\nu \over 2}} e^t}
& \mbox{for $0 < t \leqs \dd{{\nu \over 2}}$}~,
\\ 1 & \mbox{for $t > \dd{ {\nu \over 2}}$}~. \farray \right. $$
\end{prop}
\textbf{Proof.} First of all, since $v \in \Dz$ we have $e^{t \Delta} v \in \Dz$.
To go on, let us observe that
\beq \sum_{k
\in \Zd_\zz} | k |^{2 n} |(e^{t \Delta} v)_k|^2
= \sum_{k \in \Zd_\zz} | k |^{2 \nu} e^{- 2 t | k |^2} | k |^{2 n - 2 \nu} |v_k|^2 \label{thenn} \feq
$$ \leqs \Big( \sup_{k \in \Zd_\zz} | k |^{2 \nu} e^{- 2 t | k |^2} \Big)
\Big(\sum_{k \in \Zd_\zz} | k |^{2 n - 2 \nu} |v_k|^2 \Big)
\leqs \Big(\sup_{\te \in [1,+\infty)} U_{\nu t}(\te) \Big) \| v \|^2_{n-\nu}~, $$
$$ U_{\nu t}(\te):= \te^{\nu} e^{-2 t \te}~. $$
An elementary computation gives
\beq \sup_{\te \in [1,+\infty)} U_{\nu t}(\te) =
\left\{ \barray{ll} U_{\nu t}\Big(\dd{{\nu \over 2 t}}\Big) =
\dd{\left({\nu \over {2 e t}}\right)^{\nu}}  &
\mbox{for $0 < t \leqs \dd{\nu \over 2}$}~,
\\ U_{\nu t}(1) = e^{-2 t} & \mbox{for $t > \dd{\nu \over 2}$}~. \farray \right. \label{thee} \feq
Returning to \rref{thenn}, we infer that $e^{t \Delta} v \in \HM{n}$ and
\beq \| e^{t \Delta} v \|_n \leqs \sqrt{\sup_{\te \in
[1,+\infty)} U_{\nu t}(\te)} \, \| v \|_{n - \nu}~; \feq
expressing the above $\sup$ via Eq. \rref{thee} and isolating a factor $e^{-t}$ we easily get
Eqs. \rref{eqvv} \rref{eqww}. \fine
\section{Appendix. Some results on the convolution.}
\label{appeco}
The results we present here are used in the next Appendix to establish some facts about
the function $\KK_{\del}$, defined on $\Zd_\zz$ via Eq. \rref{kapde}; the cited equation
tells us that $\KK_{\del}(k)$ is a convolutionary sum (apart from
a technical detail, i.e., the elimination of $0$ and $k$ from the summation domain $\Zd$;
we return on this later on).
\parn
In this Appendix we give some general results on convolutions, which are stated independently
of the subsequent applications to $\KK_{\del}$.
Of course, the convolution of
two functions $f, g : \Zd \vain [0,+\infty)$ is defined by
\beq f * g : \Zd \vain [0,+\infty]~, \qquad k \mapsto (f * g)(k) := \sum_{h \in \Zd} f(h) g(k-h)~. \feq
We first consider the case $d=1$; our main statement for this case is contained in Proposition \ref{qcc},
which is preceded by the following definition.
\begin{prop}
\label{def1}
\textbf{Definition.} Consider a function
\beq f : \Z \vain [0,+\infty]~, \qquad k \mapsto f(k)~. \feq
$f$ is even if
\beq f(-k) = f(k) \qquad \mbox{for all $k \in \Z$}~; \feq
$f$ is unimodal if it is nondecreasing on $-\naturali = \{0,-1,-2,...\}$, and nonincreasing
on $\naturali = \{0,1,2,...\}$:
$$ f(\ell) \leqs f(m) \quad \mbox{for $\ell, m \in -\naturali$ and $\ell \leqs m$}~, $$
\beq \hspace{-0.4cm} f(\ell) \geqs f(m) \quad \mbox{for $\ell, m \in \naturali$ and $\ell \leqs m$}~. \feq
\end{prop}
Obviously enough, the unimodality of an \textsl{even} $f$ is equivalent to any one of these
conditions: $f$ is nonincreasing on $\naturali$, or
\beq f(\ell) \geqs f(m) \quad \mbox{for $\ell,m \in \Z$ and $|\ell| \leqs |m|$}~. \label{unimod} \feq
\begin{prop}
\label{qcc}
\textbf{Proposition.} Let us consider two even, unimodal functions
\beq p, q : \Z \vain [0,+\infty)~, \feq
and put
\beq s := p * q : \Z \vain [0,+\infty]~, \qquad k \mapsto s(k) = \sum_{h \in \Z} p(h) q(k-h)~. \feq
Then, $s$ is itself even and unimodal.
\end{prop}
\textbf{Proof.} (Adapted from the proof of
Proposition 4.5.5 in \cite{Lukac}; the cited result
refers to continuous convolutions,
with integrals over $\reali$ in place of sums over $\Z$). \parn
To prove that $s$ is even we write, for any $k \in \Z$:
\beq s(-k) = \sum_{h \in \Z} p(h) q(-k-h) = \sum_{j \in \Z} p(-j) q(-k+j) = \sum_{j \in \Z}
p(j) q(k-j) = s(k) \feq
(the second equality relies on the change of variable $h=-j$,
the third one holds because $p,q$ are even). \parn
Proving the unimodality of $s$ is less trivial. It suffices to prove that $s$
is decreasing on $\naturali$, which can be written as follows:
\beq s(\ell) \geqs s(\ell + 1) \qquad \mbox{for $\ell \in \naturali$}~; \label{derivs} \feq
in the sequel, we fix any $\ell$ and derive the thesis \rref{derivs}. First of all,
we write
$$  s(\ell) = \sum_{h \in \Z} p(h) q(\ell - h) = \sum_{h=-\infty}^{-1} p(h) q(\ell - h) +
\sum_{h=0}^{\infty} p(h) q(\ell - h)~; $$
now, putting $h=-k-1$ in the first sum, $h=k$ in the second one and using
$p(-k-1) = p(k+1)$, we get
\beq s(\ell) = \sum_{k=0}^{\infty} \big(~ p(k+1) q(\ell + k + 1) + p(k) q(\ell - k)~\big). \label{b8} \feq
Similarly, writing
$$  s(\ell+1) = \sum_{h \in \Z} p(h) q(\ell - h + 1) = \sum_{h=-\infty}^{0} p(h) q(\ell - h + 1) +
\sum_{h=1}^{\infty} p(h) q(\ell - h + 1)~, $$
putting $h=-k$ in the first sum, $h=k+1$ in the second one and using $p(-k) = p(k)$, we get
\beq s(\ell+1) = \sum_{k=0}^{\infty} \big(~ p(k) q(\ell + k + 1) + p(k+1) q(\ell - k)~\big). \label{b9} \feq
Now, consider any $k \in \naturali$. Then $|k| = k \leqs k + 1 = |k+1|$ and $|\ell - k | \leqs |\ell |
+ |k|$ $= \ell + k \leqs \ell + k + 1 = |\ell + k + 1|$ so that, by Eq. \rref{unimod} with $f=q$ and $f=p$,
\beq p(k) \geqs p(k+1)~, \qquad q(\ell - k) \geqs q(\ell + k + 1)~. \feq
Thus
$0 \leqs [p(k) - p(k+1)] [q(\ell - k) - q(\ell + k + 1)]$ $=
p(k) q(\ell - k) + p(k+1)$ $q(\ell + k + 1) - p(k) q(\ell + k + 1) - p(k+1) q(\ell - k)$, i.e.,
\beq p(k+1) q(\ell + k + 1) + p(k) q(\ell - k) \geqs p(k) q(\ell + k + 1) + p(k+1) q(\ell - k)~. \label{b11} \feq
Now, from \rref{b8} \rref{b9} \rref{b11} we get the thesis \rref{derivs}. \fine
Let us extend the previous considerations from the one-dimensional to the $d$-dimensional case,
for arbitrary $d$.
\begin{prop}
\textbf{Definition.} Consider a function
\beq f : \Zd \vain [0,+\infty]~, \qquad k \mapsto f(k)~; \feq
for each $r \in \{1,...,d\}$ and $k' = (k_1,...,k_{r-1},k_{r+1},...,k_d) \in \Z^{d-1}$, put
\beq f_{r k'} : \Z \vain [0,+\infty]~, \qquad k_r \mapsto f_{r k'}(k_r) :=
f(k_1,...,k_{r-1},k_r,k_{r+1},...,k_d)~. \label{effr} \feq
We call $f$ even (resp., unimodal) in each variable if, for any
$r \in \{1,...,d\}$ and $k' \in \Z^{d-1}$, the function $f_{r k'}$
is even (resp., unimodal) in the sense of Definition \ref{def1}.
\end{prop}
We note that $f$ is even in each variable if and only if
\beq f(R_r k) = f(k)~\mbox{for each $r \in \{1,...,d\}$ and $k \in \Zd$},~\feq
$$ R_r k := (k_1,...,k_{r-1},-k_r,k_{r+1},...,k_d)~. $$
If $f$ is even in each variable, recalling Eq. \rref{unimod} we see that the unimodality of $f$
in each variable is equivalent to
\beq f(\ell) \geqs f(m) \qquad \mbox{if $\ell, m \in \Zd$ and $|\ell_1| \leqs |m_1|,..., |\ell_d| \leqs |m_d|$}~.
\label{unif} \feq
\begin{prop}
\label{qccd}
\textbf{Proposition.} Let us consider two functions
\beq p, q : \Zd \vain [0,+\infty) \feq
which are even and unimodal in each variable, and put
\beq s := p * q : \Zd \vain [0,+\infty]~, \qquad k \mapsto s(k) = \sum_{h \in \Zd} p(h) q(k-h)~. \label{defes} \feq
Then, $s$ is itself even and unimodal in each variable.
\end{prop}
\textbf{Proof.} Let us fix $r \in \{1,...,d\}$, $k' = (k_1,...,k_{r-1},
k_{r+1},...,k_d) \in \Z^{d-1}$, and consider the function $s_{r k'} : \Z \vain [0,+\infty]$,
defined following Eq. \rref{effr}; we must prove that $s_{r k'}$ is even and unimodal. \parn
To this purpose we note that the definition $s := p * q$ implies
\beq s_{r k'}(k_r) = \sum_{h' \in \Z^{d-1}} (p_{r h'} * q_{r, k' - h'})(k_r) \qquad \mbox{for each
$k_r \in \Z$}~, \feq
where the functions $p_{r h'}, q_{r, k' - h'} : \Z \vain [0,+\infty)$ are defined as well following
Eq. \rref{effr}. Each of the functions $p_{r h'}, q_{r, k' - h'}$ is even and unimodal;
so, due to Proposition \ref{qcc}, their convolutions $p_{r h'} * q_{r, k' - h'}$
are even and unimodal.  A sum of even
and unimodal functions has the same properties, so our thesis about $s_{r k'}$ is proved. \fine
Our last statement about convolutions is obvious, and mentioned only for subsequent citation.
\begin{prop}
\textbf{Definition.} A function
\beq f : \Zd \vain [0,+\infty]~, \qquad k \mapsto f(k) \feq
is symmetric if, for each $\sigma$ in $\symd$ (the group of permutations of
$\{1,...,d\}$) and each $k = (k_1,...,k_d)
\in \Zd$, it is
\beq f(P_{\sigma} k) = f(k)~, \qquad P_{\sigma}(k) := (k_{\sigma(1)}, ..., k_{\sigma(d)})~. \feq
\end{prop}
\begin{prop}
\label{qcce}
\textbf{Proposition.} If $p, q : \Zd \vain [0,+\infty)$ are symmetric functions, their convolution
$p * q$ is itself symmetric.
\end{prop}
\section{Appendix. The function $\boma{\KK_{\del}}$}
\label{appeka}
Let $d \in \{2,3,...\}$. In the present
Appendix we prove a number of properties of the function $\KK_{\del}$ on $\Zd_0$, defined by Eq. \rref{kapde}; some of
these properties were mentioned in
Section \ref{sequad}.
Some of our results rely on the following two Lemmas.
\begin{prop}
\label{somme}
\textbf{Lemma.} Let us consider two radii $\rho, \rho_1$ such that $2 \sqrt{d} \leqs \rho < \rho_1 \leqs + \infty$, and
a nonincreasing function $\chi \in C([\rho - 2 \sqrt{d}, \rho_1), [0,+\infty)$). Then,
\beq \sum_{h \in \Zd, \rho \leqs | h | < \rho_1} \chi(|h|) \leqs
{2 \pi^{d/2} \over \Gamma(d/2)} \int_{\rho - 2 \sqrt{d}}^{\rho_1} dt \, (t + \sqrt{d})^{d-1} \chi(t)
\label{asinthe} \feq
(the two sides of the inequality being, possibly, $+\infty$).
\end{prop}
\textbf{Proof.} See the Appendix C of \cite{due}. \fine
\begin{prop}
\textbf{Lemma.} Fix $\nu \in (d,+\infty)$ and consider, for any $\lan \in (0,+\infty)$, the sum
\beq \Di \Ss_{\nu}(\lan) := \sum_{h \in \Zd, |h| \geqs \lan + 2 \sqrt{d}} {1 \over |h|^{\nu}}~;
\label{dsnula} \feq
this is finite, and admits the bound
\beq \Di \Ss_{\nu}(\lan) \leqs \Del \Ss_{\nu}(\lan) := {2 \pi^{d/2} \over \Gamma(d/2)}
\sum_{i=0}^{d-1} \left( \barray{c} d - 1 \\ i \farray \right) {d^{d/2-1/2-i/2} \over (\nu - i - 1)
\lan^{\nu - i - 1}}~. \label{desnu}\feq
\end{prop}
\textbf{Proof.} We apply the previous Lemma, noting that
\beq \Di \Ss_{\nu}(\lan) = \sum_{h \in \Zd, \rho \leqs | h | < \rho_1} \chi(|h|)~,
\qquad \rho := \lan + 2 \sqrt{d}, \quad \rho_1 :+ \infty, \quad \chi(t) := {1 \over t^{\nu}}~. \feq
The inequality \rref{asinthe} gives
\beq \Di \Ss_{\nu}(\lan) \leqs
{2 \pi^{d/2} \over \Gamma(d/2)} \int_{\lan}^{+\infty} dt \, {(t + \sqrt{d})^{d-1} \over t^{\nu}}~; \label{danula}\feq
now, writing the expansion $(t + \sqrt{d})^{d-1}$ $= \sum_{i=0}^{d-1}
\Big( \barray{c} d - 1 \\ i \farray \Big) \, t^i d^{d/2-1/2-i/2}$
and integrating term by term,
we see that
\beq \mbox{r.h.s.  of \rref{danula}} = \Del \Ss_{\nu}(\lan)~. \label{danula1} \feq
The relations \rref{danula} and \rref{danula1} give the thesis. \fine
From now on, we make the assumption \rref{assd}
$$ \del \in ({d \over 2} -1, +\infty)~; $$
as in Eq. \rref{kapde}, we consider the function
$$ \KK_{\del} : k \in \Zd_\zz \mapsto
\KK_{\del}(k) := \sum_{h \in \Zd_{0 k}} {1 \over |h|^{2 \del} |k - h|^2}~. $$
As noted in the comment following \rref{kapde}, the sum defining $\KK_{\del}$
certainly exists in $(0,+\infty]$; hereafter we show its finiteness,
with many other properties of the function under investigation.
\begin{prop}
\label{evaka}
\textbf{Proposition.} (i) For each $k \in \Zd_\zz$, it is $\KK_{\om}(k) < + \infty$;
furthermore, for any ``cutoff'' $\lan \in (|k|,+\infty)$, one has
\beq \Ki_{\del}(k,\lan) < \KK_{\del}(k) \leqs \Ki_{\del}(k,\lan) + \Del \Ki_{\del}(k,\lan) \label{onehas} \feq
\beq \Ki_{\del}(k,\lan) := \sum_{h \in \Zd_{0 k}, |h | < \lan + 2 \sqrt{d}}
{1 \over |h|^{2 \del} |k - h|^2}~, \label{chilan} \feq
\beq \Del \Ki_\del(k,\lan) := {2 \pi^{d/2} \over \Gamma(d/2)}
\sum_{i=0}^{d-1} \left( \barray{c} d - 1 \\ i \farray \right) {d^{d/2 - 1/2 - i/2} \over (2 \del + 1 - i)
(\lan - |k|)^{2 \del + 1 - i}}~. \label{dechilan} \feq
(ii) The function $k \mapsto \KK_{\del}(k)$ is even in each variable and symmetric, i.e.,
\beq \KK_{\del}(R_r k) = \KK_{\del}(k) \quad \mbox{for all $r \in \{1,...,d\}$ and $k \in \Zd_\zz$}~,\label{pr1} \feq
\beq \KK_{\del}(P_{\sigma}k) = \KK_{\del}(k)  \quad
\mbox{for all $\sigma \in \symd$ and $k \in \Zd_\zz$}~\label{pr2} \feq
(where, as in Appendix \ref{appeco}: $R_r k := (k_1,...,-k_r,...,k_d)$, $\symd$
are the permutations of $\{1,...,d\}$, and $P_{\sigma}k :=
(k_{\sigma(1)}, ... , k_{\sigma(d)})$). \parn
Furthermore, one has
\beq \KK_{\del}(k) + {1 \over |k|^{2}} + {1 \over |k|^{2 \del}} \leqs
\KK_{\del}(\ell) + {1 \over |\ell|^{2}} + {1 \over |\ell|^{2 \del}} \label{bob} \feq
$$ \mbox{if $k, \ell \in \Zd_{0 k}$ and $|k_1| \geqs |\ell_1|, ... , |k_d| \geqs |\ell_d|$}~. $$
(iii) It is $\sup_{k \in \Zd_\zz}\KK_{\del}(k) < + \infty$. For any $a \in \{1,2,...\}$, we have the bound
\beq \sup_{k \in \Zd_\zz}\KK_{\del}(k) \leqs \max(\SS_{\del}(a),\RR_{\del}(a))~, \label{c13} \feq
\beq \SS_{\del}(a) := \max_{k \in I(a)} \KK_{\del}(k)~, \qquad I(a) := \{ k \in \Zd_\zz~|~0 \leqs k_1 \leqs k_2 ...
\leqs k_d \leqs a \}~, \label{sia} \feq
\beq \RR_{\del}(a) :=  \KK_{\del}(0,0,...,a+1) + {1 \over (a+1)^{2}} + {1 \over (a+1)^{2 \del}}~.
\label{ria} \feq
\end{prop}
\textbf{Proof.} (i) Of course, the bounds (\ref{onehas}-\ref{dechilan}) to be proved imply finiteness
of the sum $\KK_{\del}(k)$; let us derive these bounds, for any fixed
$k \in \Zd_\zz$ and cutoff $\lan \in (|k|,+\infty)$. First of all, we write
\beq \Ki_{\del}(k,\lan) < \KK_{\del}(k) = \Ki_{\del}(k,\lan) + \Di \Ki_{\del}(k,\lan)~, \feq
where $\Ki_{\del}(k,\lan)$ is the finite sum in \rref{chilan}, and
\beq \Di \Ki_{\del}(k,\lan) :=
\sum_{h \in \Zd, |h | \geqs \lan + 2 \sqrt{d}}
{1 \over |h|^{2 \del} |k - h|^2} \label{dichilan} \feq
(note that $|h | \geqs \lan + 2 \sqrt{d}$ implies $h \neq 0, k$).
In principle, $\Di \Ki_{\del}(k,\lan) \in (0,+\infty]$; hereafter we will prove that
\beq \Di \Ki_{\del}(k,\lan) \leqs \Del \Ki_{\del}(k,\lan) \label{provt} \feq
where $\Del \Ki_{\del}(k,\lan) \in (0,+\infty)$ is defined by \rref{dechilan}; this will give
the inequality \rref{onehas} that, with the finiteness of $\Del \Ki_{\del}(k,\lan)$, also
implies $\KK_{\del}(k) < + \infty$.
In order to prove \rref{provt}, we put
\beq p := {\del + 1 \over \del}~, \qquad q := \del + 1~; \feq
then $p, q \in (1,+\infty)$ and $1/p + 1/q=1$ so that, by the H\"older inequality,
\parn
\vbox{
$$ \Di \Ki_{\del}(k,\lan) \leqs
\Big(\sum_{h \in \Zd, |h | \geqs \lan + 2 \sqrt{d}}
{1 \over |h|^{2 \del p}} \Big)^{1/p}
\Big(\sum_{h \in \Zd, |h | \geqs \lan + 2 \sqrt{d}}
{1 \over |k - h|^{2 q}} \Big)^{1/q} $$
\beq = \Big(\sum_{h \in \Zd, |h | \geqs \lan + 2 \sqrt{d}}
{1 \over |h|^{2 \del + 2}} \Big)^{1/p}
\Big(\sum_{h \in \Zd, |h | \geqs \lan + 2 \sqrt{d}}
{1 \over |k - h|^{2 \del + 2}} \Big)^{1/q}~. \label{daconsi} \feq
}
\parn
Let us consider the two sums appearing in the last passage of \rref{daconsi}. A change of
variable $h' = k - h$ in the second one gives
\beq \sum_{h \in \Zd, |h | \geqs \lan + 2 \sqrt{d}}
{1 \over |k - h|^{2 \del + 2}} = \sum_{h' \in \Zd, |h'-k | \geqs \lan + 2 \sqrt{d}}
{1 \over |h'|^{2 \del + 2}}~. \feq
On the other hand, the inequality $|h'-k | \geqs \lan + 2 \sqrt{d}$
implies $|h'| = | (h'-k) + k |$
$\geqs ||h' - k| - |k| | \geqs \lan + 2 \sqrt{d} - |k|$; so, the domain of the last
sum is contained in the domain $\{ |h'| \geqs \lan + 2 \sqrt{d} - |k| \}$, and we conclude
\beq \sum_{h \in \Zd, |h | \geqs \lan + 2 \sqrt{d}} {1 \over |k - h|^{2 \del + 2}} \leqs
\sum_{h' \in \Zd, |h' | \geqs \lan + 2 \sqrt{d} - |k|} {1 \over |h'|^{2 \del + 2}}~.
\label{ins1} \feq
Concerning the first sum in \rref{daconsi}, it is obvious that
\beq \sum_{h \in \Zd, |h | \geqs \lan + 2 \sqrt{d}}
{1 \over |h|^{2 \del + 2}}  \leqs \sum_{h \in \Zd, |h | \geqs \lan + 2 \sqrt{d} - |k|}
{1 \over |h|^{2 \del + 2}}  \label{ins2} \feq
(since the right hand side is a sum on a larger domain). We return to \rref{daconsi}, and insert therein
the bounds \rref{ins1} (with $h'$ renamed $h$) and \rref{ins2}; the conclusion is
$$ \Di \Ki_{\del}(k,\lan) \leqs \Big(\sum_{h \in \Zd, |h | \geqs \lan + 2 \sqrt{d} - |k|}
{1 \over |h|^{2 \del + 2}} \Big)^{1/p + 1/q} $$
\beq = \sum_{h \in \Zd, |h | \geqs \lan + 2 \sqrt{d} - |k|}
{1 \over |h|^{2 \del + 2}} = \Di \Ss_{2 \del + 2}(\lan - |k|) \leqs
\Del \Ss_{2 \del + 2}(\lan - |k|) ~, \label{ineds} \feq
where the last two relations follow, respectively, from the definition \rref{dsnula} of $\Di \Ss_{\nu}(\lan)$
and from the bound \rref{desnu}, here applied with
$\nu \vain 2 \del + 2$ and $\lan \vain \lan - |k|$
(note that $2 \del + 2 > d$, by the assumption \rref{assd}). On the other hand, explicitating
the definition \rref{desnu} of $\Del \Ss_{\nu}$ we see that
\beq \Del \Ss_{2 \del + 2}(\lan - |k|) = \Del  \Ki_{\del}(k,\lan)~\mbox{as in \rref{dechilan}}; \feq
with \rref{ineds}, this yields the thesis \rref{provt}. \parn
(ii) We want to show the properties (\ref{pr1}-\ref{bob}) of $\KK_{\del}$. To this purpose, let us define
\beq p : \Zd \vain [0,+\infty)~, \qquad k \mapsto p(k) :=
\left\{ \barray{ll} 1/|k|^{2 \om} & \mbox{if $k \neq 0$,} \\
1 & \mbox{if $k = 0$;} \farray \right. \feq
\beq q : \Zd \vain [0,+\infty)~, \qquad k \mapsto q(k) :=
\left\{ \barray{ll} 1/|k|^2 & \mbox{if $k \neq 0$,} \\
1 & \mbox{if $k = 0$.} \farray \right. \feq
Then, for $k \in \Zd_\zz$,
$$ (p * q)(k) =  \sum_{h \in \Zd} p(h) q(k-h) =
\sum_{h \in \Zd_{0 k}} p(h) q(k-h) + p(0) q(k) + p(k) q(0)~; $$
explicitating $p$ and $q$, we can rephrase this as
\beq (p * q)(k) = \KK_{\del}(k) +
{1 \over |k|^{2}} + {1 \over |k|^{2 \del}} \quad \mbox{for $k \in \Zd_\zz$}~. \label{eqk} \feq
The functions $p,q$ are even and unimodal in each variable, as well as symmetric
(in the sense of Appendix \ref{appeco});
by Propositions \ref{qccd}, \ref{qcce}, the same properties hold for their convolution $p * q :
\Zd \vain [0,+\infty)$. So, we have
\beq (p*q)(R_r k) = (p*q)(k) \quad \mbox{for all $r \in \{1,...,d\}$ and $k \in \Zd$}~,\label{ppr1} \feq
\beq (p*q)(P_{\sigma}k) = (p*q)(k)  \quad
\mbox{for all $\sigma \in \symd$ and $k \in \Zd$}~, \label{ppr2} \feq
\beq (p*q)(k) \leqs (p*q)(\ell)~\mbox{if $k,\ell \in \Zd$ and $|k_r| \geqs |\ell_r|$ for $r=1,...,d$} \label{bbob}
\feq
(the last relation is the inequality \rref{unif} expressing the unimodality of $f := p*q$,
with the replacement $m \vain k$).
For $k \in \Zd_\zz$, explicitating $(p*q)(k)$ via Eq. \rref{eqk} (and using the obvious relations
$|R_r k|= |P_{\sigma} k| = |k|$) Eqs. \rref{ppr1} \rref{ppr2} \rref{bbob} yield, respectively,
the conclusions  \rref{pr1} \rref{pr2} \rref{bob}. \parn
(iii) Let us prove Eq. \rref{c13} (which, of course, implies the finiteness of $\sup_{k \in \Zd_\zz}
\KK_\del(k)$). To this purpose, we choose any $a \in \{1,2,...\}$ and write \parn
{\vbox{
\beq \Zd_\zz = J(a) \cup L(a)~, \feq
$$ J(a) := \{ k \in \Zd_\zz~|~|k_r| \leqs a~\mbox{for all $r \in \{1,...,d\}$} \}~, $$
$$ L(a) := \{ k \in \Zd_\zz~|~|k_s| \geqs a+1~\mbox{for some $s \in \{1,...,d\}$}\}~. $$}
Hereafter we will derive the bounds
\beq \KK_{\del}(k) \leqs \SS_\del(a) \qquad \mbox{for all $k \in J(a)$}~, \label{c14} \feq
\beq \KK_{\del}(k) \leqs \RR_\del(a) \qquad \mbox{for all $k \in L(a)$}~, \label{c15} \feq
yielding the thesis \rref{c13}. \parn
Let $k \in J(a)$; applying a reflection to each negative component of
$k$ (if any), and then performing a suitable permutation, we can transform $k$
into an element of the set $I(a)$ in Eq. \rref{sia}; more formally, there is
a map $C = P_{\sigma} R_{r_m} ... R_{r_1}$ of $\Zd_\zz$ into itself such that
$C k \in I(a)$. Due to the results of (ii), $\KK_\del$ is invariant under $C$. So,
$$ \KK_{\del}(k) = \KK_{\del}(C k) \leqs \SS_{\del}(a)~, $$
and \rref{c14} is proved. \parn
Now, consider any $k \in L(a)$. Then $|k_s| \geqs a+1$ for some $s \in \{1,...,d\}$,
which obviously implies
$$ |k_r | \geqs |\ell_r| ~\mbox{for $r=1,...,d$}~, \qquad
\ell := (0,...,\underbrace{a+1}_{\tiny{{\mbox{place $r$}}}},...,0)~; $$
from here and \rref{bob} we infer
$$ \KK_{\del}(k) + {1 \over |k|^2} + {1 \over |k|^{2 \del}}
\leqs \KK_{\del}(\ell) + {1 \over |\ell|^2} + {1 \over |\ell|^{2 \del}}~. $$
On the other hand, the definition of $\ell$ and the symmetry of $\KK_{\del}$ give
$\KK_{\del}(\ell)$ $+ 1/|\ell|^2$ $+ 1/|\ell|^{2 \del}$
$= \KK_{\del}(0,...,0,a+1)$ $+ 1/(a+1)^2 + 1/(a+1)^{2 \del} = \RR_{\del}(a)$, so
$$  \KK_{\del}(k) +  {1 \over |k|^2} + {1 \over |k|^{2 \del}} \leqs \RR_{\del}(a)~; $$
this result is even stronger than the desired relation \rref{c15}. \fine
\section{Appendix. Evaluation of $\boma{\sup_{k \in \Z^3_\zz} \KK_{\deln}(k)}$.}
\label{appesup}
We specialize the results of the previous Appendix to the case
\beq d = 3~, \qquad \del = \delnn~; \feq
our aim is to justify the statement \rref{bounds}
$$ 27.94 < \sup_{k \in \Zd_\zz} \KK_{\deln}(k) < 32.23~. $$
First of all, let us we write down for the sup of $\KK_{\deln}$ the bound \rref{c13},
with $a=1$. In this case $I(a) = I(1) = \{(0,0,1),(0,1,1),(1,1,1)\}$, so
$\SS_{\deln}(1)$ $= \max\Big(\KK_{\deln}(0,0,1),$ $\KK_{\deln}(0,1,1),$
$\KK_{\deln}(1,1,1) \Big)$; furthermore $\RR_{\deln}(1)$
$=\KK_{\deln}(0,0,2)$ $+ 1/4 + 1/2^{7/5}$. Eq. \rref{c13} states that
$\sup_{k \in \Z^3_\zz} \KK_{\deln}$ $\leqs \max\Big(\SS_{\deln}(1), \RR_{\deln}(1)\Big)$, i.e.,
\parn
\vbox{
\beq \sup_{k \in \Z^3_\zz}\KK_{\deln}(k) \label{c19} \feq
$$ \leqs \max\Big(\KK_{\deln}(0,0,1), \KK_{\deln}(0,1,1),
\KK_{\deln}(1,1,1), \KK_{\deln}(0,0,2) + \dd{1 \over 4} + \dd{1 \over 2^{7/5}} \Big)~. $$
} \noindent
To evaluate $\KK_{\deln}$ at the points indicated above we use
Eqs. (\ref{onehas}-\ref{dechilan}), with a cutoff $\lan = 150$. The results are:
\parn
\vbox{
$$ 27.94 < \KK_{\deln}(0,0,1) < 32.23;~27.48 < \KK_{\deln}(0,1,1) < 31.77;
~26.49 < \KK_{\deln}(1,1,1) < 30.78; $$
\beq 25.69 < \KK_{\deln}(0,0,2) < 29.98; \quad \KK_{\deln}(0,0,2) + \dd{1 \over 4} + \dd{1 \over 2^{7/5}} < 30.61~.
\label{c20} \feq
}
\noindent
(Note that the reminder term $\delta \Ki_{\deln}(k,\lambda)$ of Eqs. (\ref{onehas}-\ref{dechilan})
is $O(1/\lambda^{2/5})$ for $\lambda
\vain +\infty$; such a slow decrease at infinity explains why the lower and upper bounds in \rref{c20} are not very close,
even with the fairly large chosen cutoff $\lambda = 150$.) \parn
From \rref{c19} and the upper bounds in \rref{c20} we conclude
$$ \sup_{k \in \Z^3_\zz}\KK_{\deln}(k) < 32.23~, $$
as stated in \rref{bounds}. Of course, $\KK_{\deln}(0,0,1) \leqs \sup_{k \in \Z^3_\zz}\KK_{\deln}(k)$,
and the lower bound for the former in \rref{c20} implies
$$  27.94 < \sup_{k \in \Z^3_\zz}\KK_{\deln}(k)~; $$
our justification of \rref{bounds} is concluded.
\vskip 0.4cm \noindent
\textbf{Acknowledgments.} This work was partly supported by INdAM and by MIUR, PRIN 2006
Research Project "Geometrical methods in the theory of nonlinear waves and applications".
\vfill \eject 

\noindent
\begin{thebibliography}{99}
\bibitem{Che} S.I. Chernyshenko, P. Constantin, J.C. Robinson, E.S. Titi,
\textsl{A posteriori regularity of the three-dimensional Navier-Stokes
equations from numerical computations}, J. Math. Phys. \textbf{48}(6), 065204/10 (2007).
\bibitem{DR} M. Dashti, J.C. Robinson, \textsl{An a posteriori condition on the numerical
approximation of the Navier-Stokes equations for the existence of a strong solution},
SIAM J. Numer. Anal. \textbf{46}, 3136-3150 (2008).
\bibitem{Kat} T.Kato, H. Fujita, \textsl{On the nonstationary Navier-Stokes system}, Rend. Sem. Mat.
Univ. Padova \textbf{32}, 243-260 (1962).
\bibitem{Lem} P.G. Lemari\'e-Rieusset, \textsl{Recent developments in the Navier-Stokes problem}, Chapman
\& Hall, Boca Raton (2002).
\bibitem{Lukac} E. Lukacs, \textsl{Characteristic functions}, Griffin, London (1970).
\bibitem{Pat} L. Machiels, J. Peraire, A. T. Patera,
\textsl{A posteriori finite-element output bounds for
the incompressible Navier-Stokes equations:
application to a natural convection problem},
J. Comput. Phys. \textbf{172}, 401-425 (2001).
\bibitem{uno} C. Morosi, L. Pizzocchero, \textsl{On approximate solutions of semilinear
evolution equations}, Rev. Math. Phys. \textbf{16}(3), 383-420 (2004).
\bibitem{mult} C. Morosi, L. Pizzocchero, \textsl{On the constants for multiplication in
Sobolev spaces},  Adv. in Appl. Math. \textbf{36}(4), 319-363 (2006).
\bibitem{due} C. Morosi, L. Pizzocchero, \textsl{On approximate solutions of
semilinear evolution equations II. Generalizations, and applications to Navier-Stokes
equations}, Rev. Math. Phys. \textbf{20}(6), 625-706 (2008).
\bibitem{Rob} J.C. Robinson, W. Sadowski, \textsl{Numerical verification of
regularity in the three-dimensional Navier-Stokes equations for bounded sets of
initial data}, Asymptot. Anal. \textbf{59}, 39-50 (2008).
\bibitem{Tem} R. Temam, \textsl{Navier-Stokes equations and nonlinear functional analysis}, SIAM,
Philadelphia (1983).
\bibitem{Pat2} K. Veroy, A. T. Patera, \textsl{Certified real-time solution of the parametrized steady
incompressible Navier-–Stokes equations: rigorous
reduced-basis a posteriori error bounds}, Int. J. Numer. Meth. Fluids \textbf{47}, 773-–788
(2005).
\end{thebibliography}
\end{document}